\documentstyle[12pt,amssymb]{amsart}

\newcommand\res{ \restriction }
\newcommand\exist{\exists}

\newcommand\OR{{\mbox{OR}}}

\newcommand\ult{{\mbox{ult}}}

\newcommand\mV{{\mbox{V}}}

\newcommand\Evec{{\vec E}}

\newcommand\mEvec{{{\vec {\mbox E} }}}

\newcommand\itil{{\widetilde i}}
\newcommand\jtil{{\widetilde j}}
\newcommand\ktil{{\widetilde k}}

\newcommand\gltil{{\widetilde {\gl}}}

\newcommand\Edot{{\dot E}}
\newcommand\Fdot{{\dot F}}
\newcommand\ggdot{{\dot \gg}}
\newcommand\gmdot{{\dot \gm}}
\newcommand\gndot{{\dot \gn}}

\newcommand\cMbar{{\overline {\cM}}}

\newcommand\cQbar{{\overline {\cQ}}}

\newcommand\gnbar{{\overline \gn}}
\newcommand\gzbar{{\overline \gz}}

\newcommand\cA{{\cal A}}

\newcommand\cF{{\cal F}}

\newcommand\cH{{\cal H}}

\newcommand\cJ{{\cal J}}

\newcommand\cM{{\cal M}}
\newcommand\cN{{\cal N}}

\newcommand\cP{{\cal P}}
\newcommand\cQ{{\cal Q}}

\newcommand\lb{\lbrack}
\newcommand\rb{\rbrack}

\newcommand\thra{\twoheadrightarrow}

\newcommand\lra{\longrightarrow}

\newcommand\Ra{\Rightarrow}

\newcommand\Lra{\Longrightarrow}
\newcommand\Lla{\Longleftarrow}

\newcommand\Llra{\Longleftrightarrow}

\newcommand{\LofE}{\mbox{L}\lb {\Evec} \rb}

\newcommand{\ga}{\alpha}     
\newcommand{\gb}{\beta}      
\renewcommand{\gg}{\gamma}   
\newcommand{\gd}{\delta}     
\renewcommand\ge{\varepsilon}
\newcommand{\gz}{\zeta}

\newcommand{\gk}{\kappa}  
\newcommand{\gl}{\lambda}    
\newcommand{\gm}{\mu}        
\newcommand{\gn}{\nu}

\newcommand{\gp}{\pi}        
\newcommand{\gr}{\rho}       
\newcommand{\gs}{\sigma}     
\newcommand{\gt}{\tau}

\newcommand{\go}{\omega}

\newcommand{\gG}{\Gamma}     
\newcommand{\gD}{\Delta}     
  
\newcommand{\gZ}{Z}

\newcommand{\gS}{\Sigma}

\newcommand{\ha}{\aleph}

\renewcommand{\models}{\vDash}

\newcommand{\dom}{\mbox{dom}}
\newcommand{\ran}{\mbox{ran}}
\newcommand{\crit}{\mbox{crit}}
\newcommand{\card}{\mbox{card}}

\newcommand{\cf}{\mbox{cf}}

\newcommand{\ot}{\mbox{o.t.}}

\begin{document}
\baselineskip 4ex 
\thanks{This research was partially supported by a
National Science Foundation Mathematical Sciences Postdoctoral Research
Fellowship}
\title{$\square_{\gk , < \go}$ holds in L[$\mEvec$]}{}
\author{Ernest Schimmerling}
\address{Department of Mathematics, Massachusetts Institute of Technology,
Cambridge, MA \ 02139}
\email{ernest@@math.mit.edu}

\maketitle

\ \\

\begin{center}
\underline{{\bf PRELIMINARY VERSION}}
\end{center}

\ \\

\bigskip

\def\ba{{{\mbox{{\bf a}}}}}
\def\bb{{{\mbox{{\bf b}}}}}
\def\h{{{hull}}}
\def\t{{{triv}}}
\def\l{{{lift}}}
\def\n{{{\neg\emptyset}}}
\def\e{{{\emptyset}}}

\noindent
{\bf Definition 0.}
Suppose that $\gl \leq \gk$ are cardinals 
and $\gG$ is a subset of $(\gk, \gk^+)$.  
By $\square_{\gk , < \gl} (\gG)$, we mean the principle asserting that there
is a sequence $\langle \cF_\gn \mid \gn \in \lim ( \gG ) \rangle$
such that for every $\gn \in \lim ( \gG )$,
the following hold.
\begin{list}{}{}
\item[(1)]{$1 \leq \card (\cF_\nu ) < \gl$}
\item[(2)]{The following hold for every $C \in \cF_\nu$.
	\begin{list}{}{}
	\item[(a)]{$C \subseteq \nu \cap \gG$}
	\item[(b)]{$C$ is club in $\gn$}
	\item[(c)]{$\ot (C) \leq \gk$}
	\item[(d)]{$\gnbar \in \lim (C)$ $\Ra$ $\gnbar \cap C \in
	\cF_\gnbar$}
	\end{list}}
\end{list}
By $\square_{\gk , \gl} (\gG )$ we mean $\square_{\gk , < \gl^+} (\gG )$.
If $\gG = (\gk , \gk^+ )$, 
then we write $\square_{\gk , < \gl}$ for $\square_{\gk , < \gl} (\gG)$
and 
$\square_{\gk , \gl}$ for $\square_{\gk , \gl} (\gG)$.

\bigskip

R.B.Jensen's principles $\square_\gk$ and $\square_\gk^*$ are equivalent
to $\square_{\gk , 1}$ and $\square_{\gk , \gk}$ respectively.
The strongest of these principles, $\square_\gk$, was recently
proved by Jensen to hold at all $\gk$ in the core model for non-overlapping
sequences of extenders (``below $0^{\P}$'', see \cite{J2}). 
This improves the result
in \cite{Sch}, that $\square_{\gk , 2}$ holds in the same model.
In reverse chronological order, both results build on
\cite{Wy}, \cite{We}, \cite{S},
and \cite{J1}.

Here, and in 
\cite{Sch}, the focus is on inner models of W.J.Mitchell and J.R.Steel,
introduced in
\cite{MiSt}.   These are models of the form $\LofE$, where $\vec E$ is a
sequence of partial extenders;  overlapping extenders are permitted.
The core model for one Woodin cardinal was introduced by Steel in
\cite{St1}.  This model, called K, is of the form $\LofE$.

Assume that there is no inner model with a Woodin cardinal.
That various instances of $\square_{\gk, \gl}$ hold in $\LofE$ was shown in
\cite{Sch}.  For example, $\square_\gk$ holds for all $\gk \leq \gl$,
where $\gl$ is the least measurable cardinal of Mitchell order $\gl^{++}$.
We showed that $\square_{\gk , < \gk}$ holds in $\LofE$ at every $\gk$,
and used the proof to get a lower bound of one Woodin cardinal 
on the large cardinal consistency
strength of some stationary reflection principles.
We showed that 
$\square_{\gk , cf(\gk)}$ holds for all $\gk$ in $\LofE$,
in fact, that,
$$(\gk^+)^V = (\gk^+)^{L[{\vec E}]} \Lra
V \models \square_{\gk , cf(\gk)},$$
which implies the same lower bound on the consistency
strength of the Proper Forcing Axiom.

We would like to prove that $\square_\gk$ holds in $\LofE$,
but we don't know how.
In this paper, we prove the theorem below, which subsumes the results
mentioned in the previous paragraph.
Using \cite{St2}, the theorem extends naturally to meek models 
(rather than $1$-small), and hence to the core model for a proper class
of Woodin cardinals of \cite{St3}.
We expect that the methods of \cite{SchSt} will yield extensions to tame
models, but this is not entirely worked out yet.

About the relationship between the various principles
$\square_{\gk , \gl}$, we have the following information.
Under GCH, if $\gk$ is regular, then $\square^*_\gk$ holds
(see \cite{J1}).
In \cite{BM}, it is shown that
$\square^*_{\ha_\go} \ + \ \neg \square_{\ha_\go}$ is 
consistent relative to the existence of a supercompact cardinal.
S.Todorcevic's showed that PFA
implies that $\square_{\gk , \ha_1}$ fails for 
every $\gk > \ha_1$ (the required modification of \cite{T}
was shown to the author by M.Magidor; see \cite[6.3]{Sch}).
More recently, Magidor showed that PFA is consistent with
$\forall \gk$ $\square_{\gk , \ha_2}$.
An argument due to R.M.Solovay can be adapted to show that
if $\gl$ is a limit cardinal, and $\gk \leq \gl$ is $\gl^+$-strongly compact,
then $\square_{\gl , < \gl}$ fails; moreover, if $\cf ( \gl ) < \gk$,
then $\square^*_\gk$ fails 
(this fact shown to the author by A.Kanamori; see \cite{SRK}).
Jensen showed in \cite{J2} that if there is a Mahlo cardinal,
then $\square_{\ha_1 , 2} \ + \ \neg \square_{\ha_1}$ holds 
in a forcing extension.

Corollary 2 follows from Theorem 1 using Steel's result in \cite{St1},
that the background certified core model, K$^c$, satisfies the weak
covering property ``almost everywhere''.
Corollary 3 follows directly from Theorem 1 and the weak covering
property for K proved in \cite{MiSch}.
Extensions to tame mice are possible, using the methods of \cite{SchSt},
and will appear elsewhere.

\bigskip

\noindent
{\bf Theorem 1.}
Suppose that $\LofE$ is a $1$-small weasel, all of whose 
levels are iterable.  Then for every $\gk$, $\square_{\gk , < \go}$ holds
in $\LofE$.

\bigskip

\noindent
{\bf Corollary 2.}
Assume that there is no inner model with a Woodin cardinal,
and that $\gk$ is a measurable cardinal.
Then $\square_{\gk , < \go}$ holds in V.

\bigskip

\noindent
{\bf Corollary 3.}
Assume that there is no inner model with a Woodin cardinal.
Suppose that $\gk$ is a singular cardinal such that $( \mV_{\gk^+} )^\#$ exists.
Then $\square_{\gk , < \go}$ holds in V.

\bigskip

\noindent
{\bf Proof of Theorem 1.}
We shall assume that the reader is familiar with \cite{MiSt} and 
\S\S1--4 of \cite{Sch}.  The proof here shall incorporate the ideas of \S5
of \cite{Sch}, plus a few new ones.  Unfortunately, Lemma 4.3 of
\cite{Sch} is false (an error occurs in the proof of 4.3(a)).  However,
there is a way of avoiding the use of Lemma 4.3 in the proofs of weak square
given in \S5 of \cite{Sch}.  It should be apparent from the proof that we
are about to give, what modification to \S5 of \cite{Sch} is required.

Work in $\LofE$.  Consider any cardinal $\ga$.  
Let $\gG$ be the collection of local $\ga^+$'s, that is,
the collection of limit ordinals $\gd$ such that 
$$\cJ_\gd \models \mbox{ ``$\ga$ is the largest cardinal''.}$$
For any $\gd \in \gG$, let $\gb ( \gd )$ be the least $\gb$ such that
$$\cJ_{\gb +1} \models \ \card (\gd ) = \ga.$$

For any premouse $\cM$ and for any integer $n$, by the 
$\Sigma_n$ coding structure for $\cM$, we mean the structure
$\cA_n (\cM )$ given by definition 2.7.1 or 2.7.2 of \cite{Sch}.
For any $\gd \in \gG$, we define an integer $n(\gd )$ and 
a structure $\cN (\gd )$ in
two cases, as follows.  If $\cJ_\gd$
is a type III mouse, so that $\gb ( \gd ) = \gd$, 
then we set $n(\gd ) = 1$ and 
$\cN (\gd ) = \cJ_\gd$.  
Otherwise, we let $n(\gd )$
be the least integer $n$ such that 
$\gr_n ( \cJ_{\gb (\gd )} ) = \ga$, and let $\cN (\gd )$ be the 
$\gS_{n(\gd ) -1}$ coding structure for $\cJ_{\gb ( \gd )}$.
The reason for the lack of uniformity in the definition is that
when $\cJ_\gd$ is a type III mouse, 
the $\gS_0$ coding structure for $\cJ_\gd$ is $(\cJ_\gd )^{sq}$,
the squash of $\cJ_\gd$, 
which has ordinal height $\ga$; in particular, $\gd$ is not contained
in $(\cJ_\gd )^{sq}$.  Nevertheless, in all cases $\cN (\gd )$ is an
amenable structure over which $\gd$ is $\gS_1$-collapsed.

Now fix the cardinal $\gk$ for which we are showing $\square_{\gk , <\go}$.
To avoid technicalities, assume that $\gk$ is a limit cardinal
(the other case follows already from results in \cite{Sch},
or by an easy modification of what's to come).
Let SC be the collection of $\ga < \gk$ such that $o(\ga ) \geq \gk^+$,
where $$o(\ga ) = \sup ( \{ \gn \mid \crit (E_\nu ) = \ga \})$$
(with the possibility that $o(\ga ) = \OR$ in mind).
Let $\gltil < \gk^+$ be large enough so that if $\gltil < \gn < \gk^+$
and $\crit (E_\gn ) < \gk$, then $\crit (E_\nu ) \in \mbox{SC}$.
Let $\gG = \gG_\gk - \gltil$.  Note that $\gG$ is club in $\gk^+$,
so it is enough to prove $\square_{\gk , < \go } (\gG )$, which is what
we shall do.

Given an extender $F$, let $s(F)$ be the dodd-parameter of $F$ and 
let $\gt (F)$ be the dodd-projectum of $F$.  That is, define
$s(F)$ to be the longest parameter 
$$s = \{ s_0 > \cdots > s_i \}$$ such that whenever
$0 \leq k \leq i$, we have that $s_k$ is the largest element of
$$\{ \xi > \crit (F) 
\mid \ \xi \mbox{ is an $F$-generator relative to $(s \res k)$} \},$$
and then put 
$$\gt (F)
= \sup ( \crit (F)^+ \cup \{ \xi 
\mid \xi \mbox{ is an $F$-generator relative to $s(F)$} \} ).$$
This notion is discussed in \S3 of \cite{Sch}.
In the obvious order,
$(s(F), \gt (F))$ is the least pair $(s,\gt )$ 
with 
$\gt \geq \crit (F)^+$, such that
$F$ and $F \res ( \gt \cup s)$ have the same ultrapower.
If $\cM$ is a
structure that interprets the relation symbol $\Fdot$
as coding an extender or partial extender, then we define
$s ( \cM ) = s ( \Fdot^\cM )$ and $\gt (\cM ) =  \gt (\Fdot^\cM )$.

We make some remarks and set some notation.
Recall from \cite{MiSt} that the language of potential premice (ppm's)
has only the relation symbols $\dot{\in}$, $\Edot$, and $\Fdot$,
but that the fine structure for ppm's is based on the
same language expanded by the constant symbols $\gmdot$, $\gndot$,
and $\ggdot$.  To us, a ppm will be a structure in the expanded
language, what is called an ``expanded ppm'' in \cite{MiSt}.
In fact, having the constant symbols around
only makes a significant difference in the active,
type II case.  
If $\cM$ is a passive ppm 
then $\gmdot^\cM = \gndot^\cM = \ggdot^\cM = 0$.
If $\cM$ is an active, type I ppm,
then
$\gmdot^\cM = \crit (\Fdot^\cM )$, 
$\gndot^\cM = (\gmdot^+ )^\cM$, and
$\ggdot^\cM = 0$.
And if $\cM$ is type III, then
it is $\cM^{sq}$ rather than $\cM$ over which we base the fine structure of
$\cM$.
If $\cM$ is an active ppm, then the following relationships hold.
$$\cM \mbox{ is type I } \Lra \ s( \cM ) = \emptyset \ \land
\ \gt ( \cM ) = \gndot^\cM = (\gmdot^+ )^\cM$$
$$\cM \mbox{ is type II } \Lra \  s ( \cM )_0 = \gndot^\cM - 1$$
$$\cM \mbox{ is type III }
\Lra \ s ( \cM ) = \emptyset \ \land \ \gt ( \cM ) = \gndot^\cM.$$

Partition $\gG$ into the following three pieces (no one of which need be
club in $\gk^+$).  Let $\gG_\h$ the collection of $\gn \in \gG$
for which at least one of the following hold.
\begin{list}{}{}
\item[(a)]{$n(\gn ) > 1$}
\item[(b)]{$\cJ_{\gb ( \gn )}$ is passive}
\item[(c)]{
$(\gmdot^+ )^{\cJ_{\gb ( \gn )}} \geq \gn$}
\item[(d)]{ $\cJ_{\gb ( \gn )}$ is type III and 
$\gt ( \cJ_{\gb ( \gn )} ) \geq \gn$}
\end{list}
Let $\gG_\t$ be the collection of $\gn \in \gG$ such that
$$\gn \not\in \gG_\h \ \ \land 
\ \ (\gmdot^+ )^{\cJ_{\gb ( \gn )}} < \gn \leq \gt ( \cJ_{\gb ( \gn )} ).$$
Finally, put $\gn \in \gG_\l$ if and only if $\gn \in \gG$ and
$\gt (  \cJ_{\gb ( \gn )}) = \gk$.
It is straightforward to see that $\{ \gG_\h , \gG_\t , \gG_\l \}$ is a
partition of $\gG$.
At times we shall also refer to $\{ \gG^\ga_\h , \gG^\ga_\t , \gG^\ga_\l \}$,
which we take to be the analogous partition of $\gG_\ga$.

In \S2 of \cite{Sch}, it is shown how to construct a sequence
$$B_\h = 
\langle B_\h (\gn) \mid \gn \in \lim ( \gG ) \cap \gG_\h \rangle$$ 
such that for any $\gn \in \lim ( \gG ) \cap \gG_\h$,
\begin{list}{}{}
\item[(a)]{$B_\h ( \gn )$ is club in $\gn \cap \gG$}
\item[(b)]{$\ot (B_\h (\gn ) ) \leq \gk$}
\item[(c)]{$\gnbar \in \lim (B_\h (\gn ) )$
$\Lra$ ( $\gnbar \in \lim (\gG ) \cap \gG_\h$
$\land$
$B_\h (\gnbar ) = \gnbar \cap B_\h (\gn )$ )}
\end{list}
The method for defining $B_\h$ is an extension of
Jensen's method for proving $\square_\gk$ in L, 
called the method of ``taking hulls''.
Our intuition is that
$B_\h (\gn )$ is our basic club subset 
of $\gn$ when $\gn \in \gG_\h$.
We shall define $B_\t$ shortly, and then $B_\e$ and $B_\n$ 
(both corresponding to $\gG_\l$) somewhat
later.

Using Lemma 4.4 of \cite{Sch}, we see that if $\gn \in \gG_\t$,
then 
$$\cJ_{\gb ( \gn ) + 1} \models \cf ( \gn ) = \go.$$
This means that $\square_\gk ( \gG_\t )$ is trivial.
For each $\gn \in \lim \gG_\t$, let $B_\t (\gn )$ be the
least set $C$ in the order of construction such that 
$C \subseteq \gn \cap \gG$, $C$ is cofinal in $\gn$, and
$\ot (C) = \go$.
Set 
$$B_\t = 
\langle B_\t (\gn) \mid \gn \in \lim ( \gG_\t ) \rangle.$$

We caution the reader that we are, by no means, done with the case
$\gn \in \gG_\h \cup \gG_\t$.  
We shall come back to this case again after
considering the case $\gn \in \gG_\l$.
(Remark: we shall see that $\gG_\l \subseteq \lim ( \gG )$).

Next partition $\gG_\l$ into two pieces as follows.
Let
$$\gG_\e =
\{ \gn \in \gG_\l \mid s ( \cJ_{\gb ( \gn )} ) = \emptyset \}$$
and let $\gG_\n = \gG_\l - \gG_\e$.
Then 
$$\gG_\e =
\{ \gn \in \gG_\l \mid \cJ_{\gb ( \gn )} \mbox{ is type III } \}.$$

For now, fix $\gn \in \gG_\l$
and put $\gb = \gb ( \nu )$, $F = E_\gb$,
$\gm = \crit (F)$, and
$$s = \{ s_0 > \cdots > s_\ell \} = s(F).$$
Then $\gt (F) = \gk$. Note that $\gm \in \mbox{SC}$, 
because $\gn > \gl$.
Since $\gm < \gk$, $F$ is an extender over all of
$\LofE$.  Let $W = \ult (\LofE , F )$, and let $j : \LofE \lra W$ be the
ultrapower map.

Consider an arbitrary $\gd \in \gG_\gm$.  Let 
$\cQ_\gd = \ult(\cN (\gd ), F \res ( \gk \cup s))$, and let 
$$i_\gd : \cN (\gd ) \lra \cQ_\gd$$
be the ultrapower map.
Put $\gl_\gd = (\gk^+ )^{\cQ_\gd}$.
Let $k_\gd : \cQ_\gd \lra j(\cN (\gd ))$ be given by
$$k_\gd ( [a \cup s , f]^{\cN (\gd )}_{F \res ( \gk \cup s)} )
= j(f)(a \cup s)$$
for $f \in | \cN (\gd ) |$ and $a \in [ \gk ]^{< \go}$.

\bigskip

\noindent
{\bf Remarks.}
To indicate the dependence on $\gn$ of $\cQ_\gd$, $i_\gd$, $\gl_\gd$, $k_\gd$,
etc.,
we shall sometimes write $\cQ^\gn_\gd$, $i^\gn_\gd$,
$\gl^\gn_\gd$, $k^\gn_\gd$, etc.
In \S5 of \cite{Sch}, we were writing $\gn_\gd$ for what corresponds to $\gl_\gd$
here.

\bigskip

\noindent
{\bf Lemma 1.1.}
With $\gn \in \gG_\l$,
the following hold for sufficiently large $\gd$ in $\gG_\gm$.
\begin{list}{}{}
\item[(1)]{
	\begin{list}{}{}
	\item[(a)]{$k_\gd$ is an almost 
	$\gS_1$-elementary embedding of $\cQ_\gd$
	into $$j(\cN (\gd )) = ( \cN ( j ( \gd ) ) )^W$$}
	\item[(b)]{$\gd \in \gG^\gm_\l \ 
	\Lra$\\
	$k_\gd$ is a cofinal, $\gS_1$-elementary embedding of
	$\cQ_\gd$
	into $$j(\cN (\gd )) = ( \cN ( j ( \gd ) ) )^W$$}
	\end{list}
	}
\item[(2)]{$\gk < \gl_\gd < \gn$ and $\gl_\gd = \crit ( k_\gd )$}
\item[(3)]{$\gl_\gd \in \gG$}
\item[(4)]{$\dom ( \Evec ) \cap (\gm , \gd ) \not= \emptyset$}
\item[(5)]{$0 \leq m \leq \ell \Lra
F \res ( s_m \cup (s\res m )) \in \ran (k_\delta)$}
\item[(6)]{$\cQ_\gd$ is 1-sound and 1-solid,
$\gr_1 ( \cQ_\gd ) = \gk$, and $p_1 ( \cQ_\gd ) 
= i_\gd ( p_1 ( \cN (\gd )) \cup k_\gd^{-1} (s)$}
\item[(7)]{
$\gd \not\in \gG_\e \ \Lra 
\ \cQ_\gd = \cN (\gl_\gd ) \ \land
\ n ( \gl_\gd ) = n ( \gd )$}
\item[(8)]{
\begin{list}{}{}
\item[(a)]{$\gd \in \gG^\gm_\h \Lra \gl_\gd \in \gG_\h$}
\item[(b)]{$\gd \in \gG^\gm_\t \Lra \gl_\gd \in \gG_\t$}
\item[(c)]{$\gd \in \gG^\gm_\n \Lra \gl_\gd \in \gG_\n$}
\end{list}}
\end{list}

We make a few comments before giving the proof.
With regard to (7) and (8), the case $\gd \in \gG_\e$ will
be discussed further in Lemma 1.4.
The set $F \res ( s_m \cup (s\res m ))$ in (5) is called the
$m$'th dodd-solidity witness for $F$.
In case we forget to mention it later,
clause (4) is used to know that $\cJ_{\gb (\gd )}$ satisfies the technical
hypothesis of the strong uniqueness, Theorem 6.2 of \cite{MiSt} and
of the condensation result, Theorem 2.8 of \cite{Sch}.

\bigskip

\noindent
{\bf Proof of 1.1.}
The proofs of (1)--(8) are the content of 
claims 1--5 (and environs) of the proof of Lemma 5.6 of
\cite{Sch}, but we shall remind the reader of the main ideas.

\bigskip

\noindent
{\bf Clause (1).}  Put $\gs = \sup ( j `` ( OR \cap |\cN (\gd )|))$.
Then $k_\gd$ 
is a cofinal $\gS_1$-elementary embedding 
of $\cQ_\gd$ into $j ( \cN ( \gd )) \res \gs$.
For (b), notice that $j$ is continuous at $\gb ( \gd )$
when $( \gmdot^+ )^{\cJ_{\gb ( \gd )}} < \gm$, which is the case when 
$\gd \in \gG_\l$.

\bigskip

\noindent
{\bf Clause (2).}
To assure that $\gl_\gd > \gk$,
take $\gd$ large enough so that for some $f \in J^{\Evec}_\gd$
and some $a \in [\gk ]^{<\go}$, $\gk = j(f)(a \cup s)$.
That $\gl_\gd < \gn$ follows from the amenability of $\cJ_\gb$.

\bigskip

\noindent
{\bf Clause (3).}
To be sure that $\gd \in \gG$, we take $\gd$ large enough so that
$\gl_\gd > \gltil$.  

\bigskip

\noindent
{\bf Clause (4).}
It is easy to use Theorem 8.2 of \cite{MiSt} to see that
$\dom ( \Evec ) \cap \gm^+ \not= \emptyset$.

\bigskip

\noindent
{\bf Clause (5).}
The dodd-solidity witnesses for $F$ are elements of $J^{\Evec}_\gb$,
by Theorem 3.2 of \cite{Sch}.  Take $\gd$ large enough so that 
there is a function $f \in J_\gd^{\Evec}$ and an $a \in [\gk ]^{<\go}$
such that 
$$F \res ( s_m \cup (s\res m )) 
= j(f)(a \cup s).$$

\bigskip

\noindent
{\bf Clause (6).}
First note that 
$$\cQ_\gd 
= \cH^{\cQ_\gd}_1 ( \gk \cup i_\gd (p_1( \cN (\gd ))) \cup k^{-1}_\gd (s) )$$
since $\cN (\gd )$ is 1-sound, and by the definition of $\cQ_\gd$.
So $\gr_1 ( \cQ ) = \gk$ and it 
is enough to see that the parameter 
$i_\gd (p_1( \cN (\gd ))) \cup k^{-1}_\gd (s)$ is $1$-solid over $\cQ_\gd$,
in the sense of
\cite{MiSt}.
Say that
$$p_1( \cN (\gd ) ) = p = \{ p_0 > \cdots > p_e \}$$
That $i_\gd (p)$ is $1$-solid over $\cQ_\gd$ means that
$$Th^{\cQ_\gd}_1 ( i_\gd (p_m ) \cup ( i_\gd ( p ) \res m )) \in | \cQ_\gd |$$
whenever $0 \leq m \leq e$.
But this is true by the proof of Lemma 4.6 of
\cite{MiSt}, which shows how 
$$Th^{\cQ_\gd}_1 ( i_\gd (p_m ) \cup ( i_\gd  ( p ) \res m ))$$
is $\gS_0$-definable
over $\cQ_\gd$
from 
$$\{ \mu \} \cup i_\gd  ( Th^{\cN ( \gd )}_1 ( p_m \cup ( p \res m )) );$$
note that we are using the $1$-solidity of $\cN ( \gd )$.
We must also see that
$$Th^{\cQ_\gd}_1 
( k^{-1}_\gd (s_m) \cup i_\gd ( p ) \cup (k_\gd^{-1} (s) \res m)  
\in | \cQ_\gd |$$
whenever $0 \leq m \leq \ell$.
Fix $m$ and let 
$\cP = \ult ( \cN ( \gd ), k_\gd^{-1} ( F \res ( s_m \cup (s \res m))))$,
$$\gs : \cN (\gd ) \lra \cP$$ be the ultrapower map,
and $\gp : \cP \lra j(\cN (\gd ))$ be the natural embedding
such that $j = \gp \circ \gs$.
Since $\cN (\gd )$ and $k_\gd^{-1} ( F \res ( s_m \cup (s \res m)))$
are members of $\cQ_\gd$, so is $\cP$.
A simple calculation shows that
$$Th^{\cQ_\gd}_1 ( 
k_\gd^{-1}(s_m) \cup i_\gd ( p ) \cup (k_\gd^{-1} (s) \res m))
=
Th^\cP_1 
( \gp^{-1}(s_m) \cup \gs ( p ) \cup (\gs^{-1} (s) \res m)).$$

\bigskip

\noindent
{\bf Clause (7).}
Suppose that $\gd \not\in \gG_\e$, so that
$j( \cN (\gd ) )$ is the $\gS_{n(\gd ) -1}$ coding structure for 
$ ( \cJ_{\gb ( j( \gd ) )} )^W$.
Because (1)--(6) hold, we may apply 
$\gS_n$-condensation,
Theorem 2.8 of \cite{Sch}, to conclude that
$\cQ_\gd = ( \cN (\gl_\gd ) )^W$.  But since $\LofE$ and $W$ agree
below $\gn$, $\cQ_\gd = \cN (\gl_\gd )$.

\bigskip

\noindent
{\bf Clause (8).}
By (6), we just need to examine $\cQ_\gd$.
That (a) holds is immediate from the definitions.
Suppose that $\gd \in \gG^\gm_\t$, that is, $\gd \not\in \gG_\h$
and 
$( \gmdot^+ )^{\cJ_{\gb ( \gd )}} < \gd \leq \gt ( \cJ_{\gb ( \gd )} )$.
Since $\gm = \crit (F)$ is a limit cardinal,
$$(\gmdot^+ )^{\cJ_{\gb ( \gd )}} = ( \gmdot^{\cJ_{\gb ( \gd )}} )^+
< \gm,$$ 
and so
$$(\gmdot^+ )^{\cQ_\gd} = ( \gmdot^{\cJ_{\gb ( \gd )}} )^+
< \gm < \gk < \gl_\gd.$$
By Theorem 3.2 of \cite{Sch}, $\cJ_{\gb ( \gd )}$ is dodd-amenable,
that is, for every $\xi < \gt ( \cJ_{\gb ( \gd } )$,
$$E_{\gb ( \gd )} 
\res ( \xi \cup s(E_{\gb ( \gd )} )) \in J^{\Evec}_{\gb ( \gd )}.$$
An argument essentially the same as that given in Lemma 9.1 of \cite{MiSt},
using the fact that $\gm = \crit (F) < \gt ( \cJ_{\gb ( \gd )} )$,
shows that 
$\gt ( \cQ_\gd )
= \sup ( i_\gd `` \gt ( \cJ_{\gb ( \gd )} ) )$.
Thus,
$$\gt( \cQ_\gd ) > i_\gd ( \gm ) > \gl_\gd.$$
Therefore, $\gl_\gd \in \gG_\t$.

Finally, suppose that $\gd \in \gG_\n$.
As above, 
$$(\gmdot^+ )^{\cQ_\gd} = 
( \gmdot^{\cJ_{\gb ( \gd )}} )^+
< \gm < \gk < \gl_\gd.$$
If we can show that $\gt ( \cQ_\gd ) = \gk$, then we are done.
Assuming otherwise implies that $\gt ( \cQ_\gd ) \geq \gl_\gd$,
since $\gt ( \cQ_\gd )$ is a cardinal of $\cQ_\gd$ and
$\gt ( \cQ_\gd ) \geq \gr_1 ( \cQ_\gd )$.
Then $\gl_\gd \in \gG_\t$.
As noted before, 
this implies that $\cf ( \gl_\gd ) = \go$.
But 
$$\cf ( \gl_\gd ) = \cf ( (\gmdot^+ )^{\cQ_\gd} )
= \cf ( ( \gmdot^{\cJ_{\gb ( \gd )}} )^+ ) > \go,$$
which is a contradiction.
\ \ $\blacksquare$ {\tiny \ \ Lemma 1.1}

\bigskip

\noindent
{\bf Definition 1.2.}
Let $\gD^\gn$ be the tail of $\gG_\gm$ on which the conclusion to Lemma 1.1
holds. If $\gn \in \gG_\e$, then put
$B_\e ( \gn ) = \{ \gl_\gd^\gn \mid \gd \in \gD^\gn \}$

\bigskip

When the context is clear, we write $\gD$ for $\gD^\gn$.
Notice that clauses (6)--(8) of Lemma 1.1 follow from clauses (1)--(5),
as our proof showed.  This will be important when we need to show that a
particular $\gd$ is in $\gD$.  Notice that $\gD$ is club in $\gm^+$
and that, if $\gn \in \gG_\e$, then $B_\e ( \gn )$ is club in $\gn$ of order
type $\gm^+$.

We shall need to identify $\cJ_{\gb ( \gl_\gd )}$ in the case 
$\gd \in \gG_\e$ (which was excluded in 1.1(7) and 1.1(8)).  
This is the point in \S5 of \cite{Sch} where a
reference to the fallacious
Lemma 4.3 is made, and needs to be corrected.
So suppose that $\gd \in \gG_\e$, that is, $\cJ_\gd = \cN ( \gd )$ is a 
type III mouse.  
By Lemma 1.1(1)(b), $k_\gd : \cQ_\gd \lra \cJ^W_{j(\gd )}$
is a cofinal $\gS_1$-elementary embedding.
To apply condensation,
we would need to know, at least, that $k_\gd \res |\cQ_\gd^{sq}|$
is an almost $\gS_1$-elementary  embedding of $\cQ_\gd^{sq}$
into $(\cJ^W_{j(\gd )} )^{sq}$, but this is not true.
In fact, $\cQ_\gd$ is not even a ppm, as it does not satisfy the initial
segment condition on its last extender.

\bigskip

\noindent
{\bf Definition 1.3.}
Let $\cM$ be a structure that satisfies all the conditions for being a ppm
except possibly the initial segment condition on $\Fdot^\cM$.
Put $\cP = \ult ( \cM , \Fdot^\cM )$
and let $i: \cM \lra \cP$ be the ultrapower map.
Suppose that 
$\xi < i ( \gmdot^\cM )$ and  $\xi$ is larger than any generator of
$\Fdot^\cM$.
Let $\gg = (\xi^+  )^\cP$
and let $G$ be the $( \gmdot^\cM , \gg)$-extender derived from $i$.
By the {\bf $\xi$-scaling of $\cM$},
we mean the structure
$$\langle J^\cP_\gg , \in , \Edot^\cP \res \gg , G^* \rangle$$
where 
$G^*$ is the amenable predicate coding $G$, namely
$(\gz , \gs , a , x ) \in G^*$ $\Llra$
$$(\gmdot^\cM < \gz < (\gmdot^+ )^\cM ) \ \ \&
\ \ (\xi < \gs < \gg )
\ \ \& 
\ \ (a,x) \in G \ \cap \ ([\gs ]^{<\go} \times J^\cM_\gz ) \in J^\cP_\gs$$
Note that the $\xi$-scaling of $\cM$ is not a structure in
the (expanded) language of ppm's.
By the {\bf trivial completion of $\cM$}, we mean the 
$\xi$-scaling of $\cM$ where $\xi$ is the strict supremum of the generators
of $\Fdot^\cM$.

\bigskip

\noindent
{\bf Lemma 1.4.}
Suppose that $\gn \in \gG_\l$ and $\gd \in \gG_\e \cap \gD^\gn$.
Let $\cQ^*$ be the trivial completion of $\cQ_\gd$.
Then the following hold.
\begin{list}{}{}
\item[(a)]{$\cQ_\gd$ is $1$-sound with
$s ( \cQ_\gd )  = p_1 ( \cQ_\gd ) = k_\gd^{-1} ( s )$ and
$\gt ( \cQ_\gd ) = \gr_1 ( \cQ_\gd ) = \gk$}
\item[(b)]{
$\gt ( \cQ^* ) = \gk$
and $s ( \cQ^* ) = k^{-1}_\gd ( s )$.
Moreover,
$\gr_1 ( \cQ^*) = \gk$
and
\begin{list}{}{}
\item[]{$\gn \in \gG_\e$ $\Lra$ 
$p_1 ( \cQ^*) = \emptyset$}
\item[]{
$\gn \in \gG_\n$ $\Lra$
$p_1 ( \langle \cQ^* , k^{-1}_\gd (s_0 ) \rangle ) 
=  k_\gd^{-1} ( s - \{ s_0 \} )$}
\end{list}}
\item[(c)]{
We may interpret the constant symbols $\gmdot$, $\gndot$, and $\ggdot$
in $\cQ^*$ in way that results in a premouse.  In particular,
there is an interpretation of 
$\ggdot$ in $\cQ^*$
such that
\begin{list}{}{}
\item[]{$\gn \in \gG_\e$ $\Lra$
$\langle \cQ^*, \ \gmdot^{\cQ_\gd} , \ \gk , \ 0 \ \rangle = 
\cJ_{\gb ( \gl_\gd )}$}
\item[]{
$\gn \in \gG_\n$ $\Lra$
$\langle \cQ^*, \ \gmdot^{Q_\gd} , \ k^{-1} (s_0 ) +1 , 
\ \ggdot^{\cQ^*} \ \rangle = 
\cJ_{\gb ( \gl_\gd )}$}
\end{list}
Moreover, in either case, $n ( \gl_\gd )  = 1 = n ( \gd )$
and 
$\cN( \gl_\gd ) =
\cJ_{\gb ( \gl_\gd )}$.} 
\item[(d)]{$\gl_\gd \in \gG_\e \Llra \gn \in \gG_\e$}
\item[(e)]{
$i_\gd ( \gm )$ is the least ordinal
$\xi \geq \gm$
such that
$\xi = [a , f]^{\cQ^*}_{\Fdot^{\cQ^*}}$
for some $a \in [\gm ]^{<\go}$ and 
$f : [ \gmdot^{\cQ^*} ]^{|a|} \lra \gmdot^{\cQ^*}$}
\end{list}

\bigskip

\noindent
{\bf Proof of 1.4.}

\bigskip

\noindent
{\bf Clause (a).}
From Lemma 1.1(6), we already know that 
$\cQ_\gd$ is $1$-sound with $p_1 ( \cQ_\gd ) = k_\gd^{-1} (s)$
and $\gr_1 ( \cQ_\gd ) = \gk$.
Since $k_\gd$ is cofinal and $\gS_1$-elementary, if $\xi$ is a generator of
$\Fdot^{\cQ_\gd}$, then $\xi < i_\gd ( \gm ) = k_\gd^{-1} ( j ( \gm ))$.
Thus, $\gt ( \cQ_\gd ) \leq i_\gd ( \gm )$.
Also, since $\gk$ is a cardinal, $\gt (\cQ_\gd ) \geq \gk$.

\bigskip

\noindent
{\bf Claim.}
Suppose that $\xi < i_\gd( \gm )$.
Then 
$\xi = [c \cup k^{-1}_\gd ( s ) , h ]^{\cQ_\gd}_{\Fdot^{\cQ_\gd}}$
for some 
$c \in [\gk]^{< \go}$ and 
some function
$$h : 
[\gmdot^{\cQ_\gd} ]^{|c \cup k^{-1}_\gd ( s )|} \lra \gmdot^{\cQ_\gd}.$$

\bigskip

\noindent
{\bf Proof of Claim.}
By the $1$-soundness of $\cQ_\gd$,
there is a $\gS_1$ Skolem term
$\eta$ and a parameter $b \in [\gk ]^{< \go}$
such that 
$$\xi = \eta^{\cQ_\gd} [ b, k^{-1}_\gd ( s ) ].$$
Because the ultrapower map $i_\gd$ is cofinal and $\gS_1$-elementary,
there is a $\gs < \gd$ such that 
$$\xi = \eta^{i_\gd (\cJ_\gd \res \gs )} [ b, k^{-1}_\gd ( s ) ].$$
Pick $a \in [ \gm ]^{< \go }$
and $g : \gmdot^{\cJ_\gd} \lra J^{\Evec}_{ (\gmdot^{\cJ_\gd})^+}$
such that 
$$\cJ_\gd \res \gs = [ a , g ]^{\cJ_\gd}_{\Fdot^{\cJ_\gd}};$$
this is possible because $\gd \in \gG_\e$.
Let 
$$h : [\gmdot^{\cQ_\gd} ]^{|a \cup  b \cup k^{-1}_\gd ( s ) |}
\lra \gmdot^{\cQ_\gd}$$
be defined by 
$$h(u) = \eta^{g(u^{a, a \cup  b \cup k^{-1}_\gd ( s )})} [ 
u^{b , a \cup  b \cup
k^{-1}_\gd ( s )} , u^{k^{-1}_\gd ( s ) , a \cup  b \cup k^{-1}_\gd ( s )} ].$$
Then 
$$\xi = 
i_\gd (h) (a \cup b \cup k^{-1}_\gd ( s) ) =
[ a \cup  b \cup k^{-1}_\gd ( s ) , h ]^{\cQ_\gd}_{\Fdot^{\cQ_\gd}}.$$
$\blacksquare$ {\tiny \ \ Claim}

\bigskip

It follows from the claim that $\gt (\cQ_\gd ) = \gk$ and that 
$s ( \cQ_\gd ) \leq_{lex} k^{-1}_\gd ( s )$.
But since $\cQ_\gd$ is $1$-sound,
$p_1 (\cQ_\gd )$ is the least parameter $p$ such that 
$\cQ_\gd = H^{\cQ_\gd}_1 ( \gk \cup p )$.
Recalling that
$p_1 (\cQ_\gd ) = k^{-1}_\gd ( s )$
and that 
$\cQ_\gd = H^{\cQ_\gd}_1 ( \gk \cup s ( \cQ_\gd ) )$,
we see that
$s ( \cQ_\gd )  = k^{-1}_\gd ( s )$,
and so we are done with (a).

\bigskip

\noindent
{\bf Clause (b).}
From Definition 1.3, it is clear that $s ( \cQ^* ) = s ( \cQ_\gd )
= k^{-1}_\gd ( s )$ and that $\gt ( \cQ^* ) = \gt ( \cQ_\gd ) = \gk$.
Since
$$| \cQ^* | = H_1^{\cQ^*} ( \gt ( \cQ^* ) \cup s ( \cQ^* ))
= H_1^{\cQ^*} ( \gk \cup k^{-1}_\gd  ( s )),$$
it enough to see that 
$$Th_1^{\cQ^*} ( \gk \cup q \cup \{ k^{-1}_\gd ( s_0 ) \} ) \in |\cQ^* |$$
for every parameter $q <_{lex} k^{-1}_\gd ( s - \{ s_0 \} )$.
Let $q$ be such a parameter and let $\varphi (u,v)$ be a $\gS_1$
formula in the language of  $\cQ^*$ (no constant symbols).
Then for any $c \in [ \gk ]^{< \go }$,
we have that $\cQ^* \models \varphi [ c , q ]$
if and only if
$$\exists \xi < ( k^{-1}_\gd (s_0 )^+ )^{\cQ_\gd}
\ \ \exists \eta < ( \gmdot^+ )^{\cQ_\gd}
\ \ \langle J^{\cQ_\gd}_\xi , \ \in , \ \Edot^{\cQ_\gd} \res \xi ,
\ \Fdot^{\cQ_\gd} \cap ( [\xi ]^{< \go } \times J^{\cQ_\gd}_\eta )
\ \rangle
\ \ \models
\ \ \varphi [ c , q ].$$
But
then
$Th_1^{\cQ^*} ( \gk \cup q \cup \{ k^{-1}_\gd ( s_0 ) \} )$
is $\gS_1$ definable over $\cQ_\gd$ from $q \cup \{ k^{-1}_\gd ( s_0 ) \}$.
Since $q<_{lex} k^{-1}_\gd ( s - \{ s_0 \} )$ and $p_1 ( \cQ_\gd ) =
k^{-1}_\gd (s)$,
this implies that 
$$Th_1^{\cQ^*} ( \gk \cup q \cup \{ k^{-1}_\gd ( s_0 ) \} )
\in |\cQ_\gd |.$$
But $\OR^{\cQ^*} 
= ( k^{-1}_\gd ( s_0 )^+ )^{\cQ_\gd} \geq (\gk^+)^{\cQ_\gd}$
and
$\Edot^{\cQ^*} = \Edot^{\cQ_\gd} \res \OR^{\cQ^*}$,
so we are done.

\bigskip

\noindent
{\bf Clauses (c) and (d).}
First suppose that $s = \emptyset$ (that is, $\gn \in \gG_\e$).
Then, using the elementarity of $k_\gd$ and (a),
we see that $\Fdot^{\cQ^*}$ is the trivial completion of 
$E_{j( \gd )}^W \res \gk$.
It follows from the initial segment condition on 
$\cJ^W_{j(\gd )}$ that $\cQ^*$ can be expanded to be a type III premouse,
and that the resulting premouse
is a level of $\cJ^W_{j(\gd )}$
(note that $\gk$, being a cardinal, is not the index of any extender on a
good extender sequence).
It is clear that $\cQ^*$ expands to be $\cJ^W_{\gl_\gd} = \cJ_{\gl_\gd}$.

Now suppose that $s \not= \emptyset$.
Let $\xi$ be the next generator 
of $\Edot^W_{j(\gd )}$
after $s_0$.
Note that $\xi < (s_0^+ )^W = \gb ( \gn )$, 
so that $\Edot^W \res  (\xi +1 ) = \Evec \res
(\xi +1 )$.
It follows from the initial segment condition on $\cJ^W_{j(\gd )}$
that
$E_\xi$ is the trivial completion of $E^W_{j(\gd )} \res (s_0 + 1)$.
The map $k_\gd$ restricts to a cofinal $\gS_1$-elementary embedding of
$\langle \ \cQ^* , \ \gmdot^{\cQ_\gd} , \ k^{-1}_\gd ( s_0 ) +1 \ \rangle$
into
$\cJ_\xi$
(in the language without the constant
symbol $\ggdot$).
Lemma 4.1 of \cite{Sch} implies that
$\cJ_\xi \models \eta[c , s_0 ] = \ggdot$
for some $\gS_1$ Skolem term $\eta$ (for a formula without constant symbols)
and some $c \in [\gk ]^{<\go }$.
But then
$\ggdot^{\cJ_\xi} \in \ran (k_\gd )$,
and $k_\gd$ is a cofinal $\gS_1$-elementary embedding of
$$\langle \ \cQ^* , \ \gmdot^{\cQ_\gd} , \ k^{-1}_\gd ( s_0 ) +1 , \ 
k^{-1} ( \ggdot^{\cJ_\xi} ) \ \rangle$$
into $\cJ_\xi$ 
(in the language of ppm's, including constant symbols),
and the expansion of $\cQ^*$ is a premouse.
From the definition of ppm's,
it is clear that there is a $\gS_1$ Skolem term $\gs$
(for a formula without constant symbols)
such that $\cJ_\xi \models \gs [\ggdot ] = \gndot - 1 = s_0$.
Therefore,
$$p_1 ( \langle \ \cQ^* , \ \gmdot^{\cQ_\gd} , \ k^{-1}_\gd ( s_0 ) +1 , \
k^{-1} ( \ggdot^{\cJ_\xi} ) \rangle )
= 
p_1 (\langle \ \cQ^* , \ k^{-1}_\gd ( s_0 ) \ \rangle,$$
which we saw in 1.4(b)
is equal to $k^{-1} (s - \{ s_0 \} )$.
It follows that $\cQ^*$ expands to a $1$-sound
premouse.
Applying condensation, Theorem 2.6 of \cite{Sch},
we see that $\cQ^*$ (expanded) is an initial segment of $\cJ_\xi$.
Clearly,
$\cQ^*$ (expanded) must be $\cJ_{\gb ( \gl_\gd )}$.

\bigskip

\noindent
{\bf Clause (e).}
It is enough to prove the statement with $\cQ^*$ replaced by $\cQ_\gd$.
But then the statement is clear, since
$$\{ [a , f ]^{\cQ_\gd}_{\Fdot^{\cQ_\gd}} \mid a \in [\gm ]^{< \go}
\ \land
\ f: [ \gmdot^{\cQ_\gd} ]^{|a|} \lra \gmdot^{\cQ_\gd} \}
= \ran ( i_\gd ) \cap \OR^{\cQ_\gd}.$$
$\blacksquare$ {\tiny \ \ Lemma 1.4}

\bigskip

We are still considering an arbitrary $\gn \in \gG_\l$,
with $\gb = \gb ( \gn )$, $F = E_\gb$, $\gm = \crit (F)$, etc.
(there is no assumption on $\gd$).

\bigskip

\noindent
{\bf Definition 1.5.} Suppose that $E$ is an extender over $\cM$, 
with $\gm = \crit (E)$ and $s = s(E)$. 
Let $\gd \in (\gm , ( \gm^+ )^\cM )$.  
We say that $\gd$ is {\bf $E$-closed} 
if and only if for every $x \subseteq \gm$, 
if there are $a \in [\gd ]^{<\go}$ and $f \in J^\cM_\gd$ such that
$x = [a \cup s , f ]^\cM_E$, then $x \in J^\cM_\gd$.

\bigskip

\noindent
{\bf Lemma 1.6.}
Suppose that $\gn \in \gG_\l$, $\gd \in \gD$, and $\gd$ is $F$-closed.
Let $r$ be the longest initial segment $q$ of 
$p_1 ( \cQ_\gd )$
such that 
if $\cM = \cH^{\cQ_\gd}_1 ( \gm \cup q )$
and $\gp : \cM \lra \cQ_\gd$ is the inverse of the transitive collapse,
then
$\gp^{-1}(q ) = p_1 ( \cM )$.
Then $r = i_\gd ( p_1 ( \cN ( \gd )))$.

\bigskip

\noindent
{\bf Proof of 1.6.}
Recall that
$p_1 ( \cQ_\gd ) = i_\gd ( p_1(\cN ( \gd ))) \cup k^{-1}_\gd ( s )$.
So $r$ as defined in the statement of the lemma is at least as long at 
$i_\gd ( p_1(\cN ( \gd )))$.

\bigskip

\noindent
{\bf Claim.}
If $\cH = \cH_1^{\cQ_\gd} ( \gm \cup p_1 ( \cQ_\gd ))$,
then $( \gm^+ )^\cH = \gd$.

\bigskip

\noindent
{\bf Proof of Claim.}
Suppose that $x \subseteq \gm$ and $x \in | \cH |$.
Then there is a
$\gS_1$-Skolem term $\eta$, a parameter $c \in [\gm ]^{< \go}$,
and an ordinal 
$\gs < \OR^{\cN ( \gd )}$ such that
$$x = \gm \cap \eta^{i_\gd ( \cN ( \gd ) \res \gs )}
[ c , i_\gd ( p_1 ( \cN ( \gd ))) , k^{-1}_\gd (s) ].$$
This is by the $1$-soundness of $\cN (\gd )$ and the fact that $i_\gd$ is a
cofinal, $\gS_1$-elementary embedding.
Let 
$$f(u) = \eta^{\cN ( \gd ) \res \gs} 
[c , p_1 ( \cN (\gd )), (u - u_0) ] \cap u_0$$
defined on $[\gm ]^{|s| +1}$.
Then $f \in J^{\vec E}_\gd$ and $x = i_\gd (f)(s \cup \{\gm \})
= [s \cup \{\gm \} , f]^{J^\Evec_\gd}_F$.
Since $\gd$ is $F$-closed, $x \in J^\Evec_\gd$
\ \ $\blacksquare$ {\tiny \ \ Claim.}

\bigskip

For contradiction, assume that $r$
properly extends $i_\gd ( p_1 ( \cN ( \gd )))$.
Let $$\cM = \cH^{\cQ_\gd}_1 ( \gm \cup r )$$ and let 
$\gp: \cM \lra \cQ_\gd$ be the inverse of the transitive collapse.
By the claim and the fact that $r$ contains $i_\gd (p_1 ( \cN ( \gd )))$,
we have that
$$\gd = (\gm^+ )^{\cN (\gd )} \leq (\gm^+ )^\cM \leq (\gm^+ )^\cH 
= \gd$$
(where $\cH$ is as in the claim).
So $(\gm^+ )^\cM = \gd$.
By condensation, $\cM$ is an initial segment of $\cQ_\gd$ and hence of
$\LofE$.  But then
$\cM = \cN ( \gd )$.  
This is absurd, since
$p_1 ( \cM ) = \pi^{-1} ( r )$ is strictly longer than $p_1 (\cN ( \gd ))$.
\ \ $\blacksquare$ {\tiny \ \ Lemma 1.6}

\bigskip

\noindent
{\bf Definition 1.7.}
Let $\cM$ be a $1$-sound, amenable structure.
Suppose that $ \ga < \gr_1 ( \cM )$
Let $r_\ga (\cM )$
be the longest initial segment $q$ of $p_1 (\cM )$ such that 
if $\gp : \cH \lra \cM$
is the inverse of the transitive collapse for $H^\cM_1 ( \ga \cup q )$,
then 
$\gp^{-1} ( q ) = p_1 (\cH )$.
We say that $\ga$ is {\bf switching for $\cM$} if
$r_\gb ( \cM )$ 
is strictly longer than $r_\ga ( \cM )$ whenever $\ga < \gb$.

\bigskip

\noindent
{\bf Lemma 1.8.}
With $\gn \in \gG_\l$, $\gd \in \gD$ and $\gd$ $F$-closed,
we have that
$$r_\gm ( \cQ_\gd ) = i_\gd ( p_1 ( \cN (\gd ))).$$
Moreover, if $\gn \in \gG_\n$ and $\gd$ is sufficiently large, then $\gm$
is switching for $\cQ_\gd$.

\bigskip

\noindent
{\bf Proof of 1.8.}
That $r_\gm ( \cQ_\gd ) = i_\gd ( p_1 ( \cN (\gd )))$
is just a restatement of Lemma 1.6.  Assume that $\gn \in \gG_\n$,
so that $s \not= \emptyset$.  Let 
$$X = 
H^{\cQ_\gd}_1 ( \gm^+ \cup i_\gd (p_1 ( \cN (\gd ))) \cup \{ s_0 \} )$$
(the uncollapsed hull).

\bigskip

\noindent
{\bf Claim.}
If $\gd$ is sufficiently large,
then
$F \res s_0 \in \ran ( k_\gd )$ and 
$k^{-1}_\gd ( F \res s_0 ) \in X$.

\bigskip

\noindent
{\bf Proof of Claim.}
The initial segment condition on $\cJ_\gb$ implies that 
$F \res s_0 \in J^{\Evec}_\gb$, so for large $\gd$ we have that
$F \res s_0 \in \ran (k_\gd )$.
Applying
Lemma 4.1 of \cite{Sch} to $\cJ_\gb$ gives a 
$\gS_1$-Skolem term $\eta$ and a parameter
$a \in [\gm^+ ]^{< \go}$ such that
$$F \res s_0 = \eta^{\cJ_\gb} [ \gm  , a , s_0].$$
It is important here that $\eta$ is a Skolem term for a formula 
without constant symbols.
Fix some $\gg < \gb$ such that
$F \res s_0 = \eta^{\cJ_\gb \res \gg } [ \gm  , a , s_0].$
Assume that $\gd$ is large enough so that
$\cJ_\gb \res \gg \in \ran ( k_\gd )$.
It is enough to see $k_\gd^{-1} ( \cJ_\gb \res \gg ) \in X$,
for then
$$k^{-1}_\gd ( F \res s_0 ) = 
\eta^{k_\gd^{-1} ( \cJ_\gb \res \gg ) } [ \gm  , a , k^{-1}_\gd (s_0 )]
\in X$$
The last predicate of $\cJ_\gb \res \gg$ codes an extender fragment,
so there is some $\gd_0 < \gm^+$ such that
$\Fdot^{\cJ_\gb \res \gg }$
codes
$$G = F \cap ( [\gg ]^{<\go } \times J^\Evec_{\gd_0}).$$
Without loss of generality,
$\gd_0 \in \gD$ and 
$$\ult ( J^\Evec_{\gd_0} , G ) \models s_0^+ = \gg.$$
Then $G \in \ran (k_\gd )$ and it is enough to see that
$k^{-1}_\gd  (G \res (s_0 +1)) \in X$.
Now $G \res (s_0 +1)$ is the $(\gm , s_0 + 1 )$-extender derived from
$j \res |\cN ( \gd_0 )|$, which in turn is
the inverse transitive collapsing map:
$$\cN (\gd_0 ) = 
\cH^{j ( \cN (\gd_0 ))}_1 ( \gm \cup j(p_1 (\cN ( \gd_0 )))) \lra 
j (\cN (\gd_0 )).$$
If we assume that $\gd > \gd_0$, then
$k^{-1}_\gd (G \res (s_0 +1))$ 
is the $(\gm , k^{-1}_\gd (s_0)+1)$-extender derived
from the 
inverse transitive collapsing map:
$$\cN (\gd_0 ) =
\cH^{i_\gd ( \cN (\gd_0 ))}_1 ( \gm \cup i_\gd (p_1 (\cN ( \gd_0 )))) \lra 
i_\gd (\cN (\gd_0 )).$$
Since $\ran (i_\gd ) \cup \{ \gm , k^{-1}_\gd  (s_0) \} \subseteq X$, we are done.
\ \ $\blacksquare$ {\tiny \ \ Claim}

\bigskip

Thus both $k_\gd^{-1} ( F \res s_0 )$ and $\cN (\gd )$ are elements of 
$X$.  By the computation done to prove Lemma 1.1(6), we have that
$$Th^{\cQ_\gd}_1 ( k^{-1}_\gd (s_0 ) \cup i_\gd ( p_1 ( \cN (\gd ))) ) \in X.$$
That is, the $(\ell +1)$'st solidity witness for $\cQ_\gd$ is a member of
$X$.
Therefore,
$$i_\gd ( p_1 ( \cN (\gd ))) \cup \{ k^{-1}_\gd (s_0) \}$$
is an initial segment of $r_{\gm^+} ( \cQ_\gd )$,
and the lemma follows.
\ \ $\blacksquare$ {\tiny \ \ Lemma 1.8}

\bigskip

\noindent
{\bf Definition 1.9.}
If $\gn \in \gG_\n$, then let 
$\Phi^\gn$
be the tail of $\gd$ in
$\{ \gd \in \gD^\gn \mid \gd \mbox{ is $F$-closed} \}$ 
for which
$\gm$ is switching for $\cQ_\gd$.
Put
$B_\n ( \gn ) = \{ \gl_\gd^\gn \mid \gd \in \Phi^\gn \}$.

\bigskip

When the context is clear, we write $\Phi$ for $\Phi^\gn$.
Notice that if $\gn \in \gG_\n$, then $\Phi$ is club in $\gm^+$ and 
$B_\n (\gn )$ is club in $\gn$ of order type $\gm^+$.
We are done defining the basic sequences $B_\h$, $B_\t$,
$B_\e$, and $B_\n$.  Next, we go back to considering an arbitrary $\gn \in
\gG$.

\bigskip

\noindent
{\bf Definition 1.10.}
For each $\gn \in \gG$, we define a finite set,
$\cM ( \gn )$, of amenable, $1$-sound structures.
If $\gn \in \gG - \gG_\l$, 
then let $\cM ( \gn ) = \{ \cN (\gn ) \}$.
Suppose now that $\gn \in \gG_\l$.  
For each $\ga < \gk$, define $\gg_\ga$ to be the least ordinal $\gg$
such that
$\gg = [a , f ]^{\cJ_{ \gb ( \gn )} }_{ \Fdot^{\cJ_{\gb ( \gn )}} }$
for some 
$a \in [\ga ]^{< \go}$
and
$f: [ \gmdot^{\cJ_{\gb ( \gn )}} ]^{|a|} \lra \gmdot^{\cJ_{\gb ( \gn )}}$.
Notice that $\ga \mapsto \gg_\ga$ is non-increasing, and therefore
takes on only finitely many values. 
Whenever 
$ \gmdot^{\cJ_{\gb ( \gn )}} < \ga < \gk$, let
$\cM_\ga$ be the $\gg_\ga$-scaling of
$\cJ_{\gb ( \gn )}$.
Let $\cM ( \gn )$ be the finite set
$$\{ \cJ_{\gb ( \gn )} \} \cup
\{ \cM_\ga \mid  
\gmdot^{\cJ_{\gb ( \gn )}} < \ga \in \mbox{SC}
\ \land \ \cM_\ga \mbox{ is $1$-sound}
\ \land \ p_1 (\cM_\ga ) = s (E_{\gb ( \gn )}) \}.$$

\bigskip

The idea behind Definition 1.10 is as follows. 
Suppose that $\gn \in \gG$ and that $\cM \in \cM (\gn )$.
Using $\cM$, we shall define another collection of club subsets of $\gn$
that are intended to anticipate the possibility that $\gn$ is 
a limit point of the basic club subset for some $\gn' \in \gG_\l - \gn$.
When $\gn \in \gG - \gG_\l$, $\cM$ is just $\cN ( \gn )$.
When $\gn \in \gG_\l$,
we are allowing for the possibility that
there is some $\gn' \in \gG_\l -  \gn$ with 
$\crit (E_{\gb ( \gn' )}) = \ga$, and some
$\gd \in \gG^\ga_\e$, such that $\gl^{\gn^{\prime}}_\gd = \gn$ and
$i_\gd^{\gn'} ( \ga ) = \gg^\gn_\ga$,
and
$\cQ^{\gn'}_\gd = \cM^\gn_\ga \in \cM ( \gn )$.
Since we are supposed to be proving $\square_{\gk , < \go}$,
it is important that $\cM (\gn )$ be finite.
The next lemma is immediate from Lemmas 1.1 and 1.4, and should make the
point clearer.

\bigskip

\noindent
{\bf Lemma 1.11.}
If $\gn \in \gG_\l$ and $\gd \in \gD^\gn$,
then 
$\cQ^\gn_\gd \in
\cM (\gl^\gn_\gd )$

\bigskip

\noindent
{\bf Definition 1.12 (long).}
Suppose that $\gn \in \gG$ and $\cM \in \cM (\gn )$.
Suppose that $\ga < \gk$
and $q$ is an initial segment of
$p_1 ( \cM)$.  
Put $s_{\ga , q}  = p_1(\cM) - q$.
Define $\cH^\cM_{\ga , q} = \cH^\cM_1 ( \ga \cup q )$ and
$j^\cM_{\ga , q} : \cH^\cM_{\ga , q} \lra \cM$ 
to be the inverse of the transitive collapse,
and $\gd^\cM_{\ga , q} = ( \ga^+ )^{\cH^\cM_{\ga , q}}$.
We say that $\ga$ is {\bf $\cM$-$q$-reasonable} if and only if
\begin{list}{}{}
\item[(i)]{$\ga$ is a limit cardinal}
\item[(ii)]{$\crit (j^\cM_{\ga , q} ) = \ga$ 
and $j^\cM_{\ga , q} ( \ga ) > \max ( s \cup \{\gk \} ) $}
\item[(iii)]{$j^\cM_{\ga , q}$ is cofinal in the ordinals of $\cM$}
\item[(iv)]{$\gG_\ga \cap \gd^\cM_{\ga , q}$ is cofinal in $\gd^\cM_{\ga , q}$}
\item[(v)]{$\dom (\Evec ) \cap (\ga , \gd^\cM_{\ga , q} ) \not= \emptyset$}
\end{list}
Assume that at least condition (ii) holds.
Define $F^\cM_{\ga , q}$ to be the superstrong extender fragment
derived from $j^\cM_{\ga , q}$, that is,
$$F^\cM_{\ga , q} 
= \{ (a, x) \in [j^\cM_{\ga , q} ( \ga ) ]^{<\go} \times |\cH^\cM_{\ga , q}|
\mid 
x \subset [\ga ]^{|a|} \ \land \ a \in j^\cM_{\ga , q} (x) \}.$$
Given any $\gd \in \gG_\ga \cap \gd_{\ga , q}$,
we define 
$$( \cQ_\gd^{\ga , q} )^\cM 
= \ult (\cN (\gd ) , F^\cM_{\ga , q} \res ( \gk \cup s ))$$
and
$$( i^{\ga , q}_\gd )^\cM : \cN (\gd ) \lra ( \cQ_\gd^{\ga , q} )^\cM$$
to be the ultrapower map
and
$$( k_\gd^{\ga , q} )^\cM : ( \cQ_\gd^{\ga , q} )^\cM \lra 
j^\cM_{\ga , q} ( \cN (\gd ) )$$
to be the natural map such that
$$j^\cM_{\ga , q} = ( k_\gd^{\ga , q} )^\cM \circ ( i^{\ga , q}_\gd )^\cM .$$
Define
$$( \gl_\gd^{\ga , q} )^\cM = ( \gk^+ )^{ ( \cQ_\gd^{\ga , q} )^\cM }.$$
Finally (for Definition 1.12),
we say that
a pair $(\ga , q)$ is {\bf $\cM$-considered} if and only if either
$$\mbox{$\ga$ is $\cM$-$q$-reasonable,}$$ 
or
$$\exist \gn \in \gG_\l \ \ \cM = \cJ_{\gb ( \gn )} \mbox{ and }
(\ga , q) = ( \gmdot^\cM , \emptyset ).$$

\bigskip

\noindent
{\bf Notation.}
Note that if $\gn \in \gG_\l$, then $\gmdot^{\cJ_{\gb ( \gn )}} < \ga$
whenever $\ga$ and $q$ satisfy condition 1.12(ii),
because $\gmdot^{\cJ_{\gb ( \gn )}}$ is in every hull.
We shall only be interested in $\ga$ that satisfy at least condition
1.12(ii), and so we can make the
notational change below without fear of
ambiguity with Definition 1.12.
Suppose that $\gn \in \gG_\l$.  Then:
\begin{tabbing}
xxxxxxxxxxxxxxx\= xxxxxxxxxxxxxxxxx\= xxxxxxxxxxxxxx\= \kill
\> \underline{{\sc old notation}}: \> \underline{{\sc new notation}}: \\
\> \ \ $F^\gn$		\> \ \ $F^{\cJ_{\gb ( \gn )}}_{\gm , \emptyset}$\\
\> \ \ $j^\gn$		\> \ \ $j^{\cJ_{\gb ( \gn )}}$\\
\> \ \ $s^\gn$ 		\> \ \ $s^{\cJ_{\gb ( \gn )}}_{\gm , \emptyset}$\\
\> \ \ $(\gmdot^+)^{\cJ_{\gb ( \gn )}}$
			\> \ \ $\gd^{\cJ_{\gb (\gn )}}_{\gm ,\emptyset}$\\
\> \ \ $\cQ^\gn_\gd$
		\> \ \ $(\cQ_\gd^{\gm ,\emptyset})^{\cJ_{\gb ( \gn)}}$\\
\> \ \ $i^\gn_\gd$      
		\> \ \ $(i^{\ga ,\emptyset}_\gd )^{\cJ_{\gb (\gn )}}$\\
\> \ \ $\gl^\gn_\gd$        
		\> \ \ $(\gl^{\ga ,\emptyset}_\gd )^{\cJ_{\gb (\gn )}}$\\
\> \ \ $k^\gn_\gd$      
		\> \ \ $(k_\gd^{\gm ,\emptyset})^{\cJ_{\gb ( \gn)}}$\\
\> \ \ $\gD^\gn$	\> \ \ $\gD^{\cJ_{\gb (\gn )}}_{\gm ,\emptyset}$\\
\> \ \ $\Phi^\gn$	\> \ \ $\Phi^{\cJ_{\gb (\gn )}}_{\gm ,\emptyset}$\\
\end{tabbing}

As we often dropped the superscript ``$\gn$'' in the old notation,
we shall drop the superscript ``$\cJ_{\gb ( \gn )}$'' or ``$\cM$''
in the new notation
when the context is clear.
For us, the two most important instances of
$\ga$ being $q$-$\cM$-reasonable will
be
$q = p_1 (\cM )$ and
and $q = r_\ga ( \cM )$.

The definitions of the dodd-parameter and dodd-projectum extend
in an obvious way to extender fragments (see \cite{MiSchSt}).
The following lemma implies that in all cases,
$s(F_{\ga , q}) = s_{\ga , q}$ and 
$\gt (F_{\ga , q}) = \gk$, and that $F_{\ga , q}$ is dodd-solid,
dodd-amenable, and weakly-amenable.

\bigskip

\noindent
{\bf Lemma 1.13.}
Suppose that
$\gn \in \gG$, $\cM \in \cM (\gn )$, and $\ga$ is $\cM$-$q$-reasonable.
Let $s = s_{\ga , q} = p_1 (\cM ) - q$, $j = j_{\ga , q}$, 
$\cH = \cH_{\ga , q}$,
and $F = F_{\ga , q}$.
Say $s = \{ s_0 > \cdots > s_\ell \}$.
Then the following hold.
\begin{list}{}{}
\item[(a)]{$\ult ( \cH , F \res ( \gk \cup s )) = 
\ult ( \cH , F ) = \cM$}
\item[(b)]{$0 \leq m \leq \ell 
\Lra F \res (s_m \cup (s \res m )) \in |\cM |$}
\item[(c)]{$\forall \xi < \gk$ \ \ $F \res ( \xi \cup s ) \in |\cM |$}
\item[(d)]{$\forall \gd < \gd_{\ga , q}
\ \ F \cap ( [\gk \cup s ]^{< \go } \times J^\cM_\gd ) \in |\cM |$}
\end{list}

\bigskip

\noindent
{\bf Proof of 1.13.}

\noindent
{\bf Clause (a).}
Suppose that $x \in |\cM |$.
Since $\cM$ is $1$-sound, there is a $\gS_1$-Skolem term
$\eta$, and a parameter $a \in [\gk ]^{< \go }$ such that
$$x = \eta^\cM [ a , s , q ].$$
By 1.12(iii), there is a $\gs < \OR^\cH$ such that 
$$x = \eta^{\cM \res j ( \gs )} [ a , s , q ].$$
Define $f : [ \ga ]^{|a \cup s |} \lra |\cH | $ by
$$f ( u ) = \eta^{\cH \res \gs} [ u^{a, a\cup s} , u^{s, a \cup s} , j^{-1}
(q )].$$
Then $f \in |\cH |$ and $x = [a \cup s , f ]^\cH_F = j(f)(a \cup s)$.

\bigskip

\noindent
{\bf Clauses (b) and (c).}
Suppose that $0 \leq m \leq \ell$.
Since $\cM$ is $1$-solid,
$$Th^\cM_1 ( s_m \cup q \cup (s \res m )) \in |\cM |.$$
From this theory, we can easily compute 
$F \res (s_m \cup (s \res m ))$ inside $\cM$.
Clause (c) follows from the fact that
$\gp_1 (\cM ) = \gk$.

\bigskip

\noindent
{\bf Clause (d).}
Given $\gd < \gd_{\ga , q}$, pick $e \in |\cH |$ 
such that $e: \ga \thra \cP ( \ga ) \cap \cH_{\ga , q}$.
Then for $a \in [\gk ]^{<\go}$ and $\xi < \ga $,
$$(a \cup s , e(\xi ) ) \in F \ \Llra \ a \cup s \in j(e)(\xi)$$
Thus
$F \cap ( [\gk \cup s ]^{< \go } \times J^\cM_\gd ) \in |\cM |$
\ \ $\blacksquare$ {\tiny \ \ Lemma 1.13}

\bigskip

The following lemma is proved essentially as were Lemmas 1.1 and 1.4,
using Lemma 1.13 in the obvious way.

\bigskip

\noindent
{\bf Lemma 1.14.}
Suppose that
$\gn \in \gG$, $\cM \in \cM (\gn )$, and $(\ga , q)$ is $\cM$-considered.
Then, for sufficiently large $\gd \in (\ga , \gd_{\ga , q}^\cM )$,
all clauses of Lemmas 1.1 and 1.4 hold for $(\ga , q)$,
in the obvious sense.
In particular,
$\cQ_\gd^{\ga , q} \in \cM ( \gl^{\ga , q}_\gd )$.

\bigskip

We next define a new family of club subsets of $\gn$,
one for each triple
$(\cM , \ga , q )$ such that $\cM \in \cM (\gn )$,
$q$ an initial segment of $p_1 (\cM )$,
and $\ga$ is $\cM$-$q$-reasonable.
It seems likely that for a fixed $\cM$, 
there could be infinitely many corresponding
$\ga$;
this entails a difficulty to be dealt with later.

\bigskip

\noindent
{\bf Definition 1.15.}
Suppose that $\gn \in \gG$, $\cM \in \cM (\gn )$ and
$(\ga , q)$ is $\cM$-considered.
Define $\gD_{\ga , q}^\cM$ to be the tail of 
$\gG_\ga \cap \gd^\cM_{\ga , q}$ on
which the conclusion of Lemma 1.14 holds.
Define
$$D (\cM , \ga , q)
=
\{ ( \gl^{\ga , q}_\gd )^\cM 
\mid \gd \in \gD_{\ga , q}^\cM \}$$

\bigskip

Notice that $\cM$ determines $\gn$.
Also notice that 
Definition 1.15 
is consistent with the definition of $\gD^\cM_{\gm , \emptyset}$
given in the case $\gn \in \gG_\l$ and $\cM = \cJ_{\gb ( \gn )}$.
If $\gn \in \gG_\l$ with $\gm = \crit (E_{\gb ( \gn)})$,
then, by definition,
$$\gn \in \gG_\e \Lra 
B_\e (\gn )= D(\gn , \cJ_{\gb ( \gn )} , \gm , \emptyset )$$
$$\gn \in \gG_\n \Lra
B_\n (\gn ) \subseteq D(\gn , \cJ_{\gb ( \gn )} , \gm , \emptyset ).$$
The following lemma follows easily from what we have already established.
Together with the existence of our sequences $B_\h$ and $B_\t$,
Lemma 1.16 already implies that
$\square_{\gk , \gk} ( \gG )$ holds.

\bigskip

\noindent
{\bf Lemma 1.16.}
Suppose that $\gn \in \gG$ and $\cM \in \cM (\gn )$.  
Assume that $(\ga , q)$ is $\cM$-considered.
Then the following hold.
\begin{list}{}{}
\item[(a)]{
$\gd \mapsto \gl^{\ga , q}_\gd$ is a continuous,
increasing, cofinal function
from $\gD_{\ga , q}$ into $\gn$.
In particular, 
$D (\cM , \ga , q )$ is a club subset of $\gn$
with
$$ \ot ( D (\cM , \ga , q ) ) \leq \gd_{\ga , q} \leq \ga^+ \leq \gk.$$}
\item[(b)]{
Suppose that $\gd \in \lim ( \gD_{\ga , q} )$. 
Let $p = i^{\ga , q}_\gd (p_1 (\cN ( \gd )))$.
Then 
$\ga$ is 
$\cQ_\gd^{\ga , q}$--$p$--reasonable
and 
$$D (\cQ^{\ga , q}_\gd , \ga , p )
= \gl^{\ga , q}_\gd \cap
D (\cM , \ga , q ).$$}
\end{list}

\bigskip

In Lemmas and Definitions 1.17--1.21, we are mainly concerned with the case 
in which $s_{\ga , q} = \emptyset.$

\bigskip

\noindent
{\bf Lemma 1.17.}
Suppose that $\gn \in \gG$ and $\cM \in \cM (\gn )$.  
Assume that $\ga$ is $\cM$-$q$-reasonable,
where $q$ is an initial segment of $p_1 (M )$,
and $s = s_{\ga ,q} = p_1 (\cM ) - q$.
Then the following hold.
\def\bg{{{\mbox{{\bf g}}}}}
\begin{list}{}{}
\item[(a)]{Suppose that $\gb > \ga$ and $\gb$ is $\cM$-$q$-reasonable.
Suppose further that 
$j_{\gb , q} ( \gb ) = j_{\ga , q} ( \ga )$.
Then 
$$D (\cM , \ga , q )
\subseteq
D (\cM , \gb , q ).$$ }
\item[(b)]{
The collection of $\gb$ such that
$\gb$ is $\cM$-$q$-reasonable and
$j_{\gb , q} ( \gb ) = j_{\ga , q} ( \ga )$}
\item[(c)]{
If $q = p_1 (\cM)$, then there is a largest $\gb$ such that $\gb$ is 
$\cM$-$q$-reasonable.}
\item[(d)]{
Suppose that $q = p_1 ( \cM )$.
Suppose that $\gd \in \gD_{\ga , q}$ and $\gb \in (\ga ,\gk )$
are such that $\gk = j_{\ga , q} (f ) (a)$
for some $a \in [\gb ]^{<\go}$ and some 
$f: [\ga ]^{|a|} \lra \ga$.
Let $p = i^{\ga , q}_\gd ( p_1 (\cN (\gd )))$.
Then
$\ga$ is $\cQ^{\ga , q}_\gd$--$p$-reasonable,
but
$\gb$ is not $\cQ^{\ga , q}_\gd$--$p$-reasonable.}
\end{list}

\bigskip

\noindent
{\bf Proof of Lemma 1.17.}

\noindent
{\bf Clause (a).}
Fix an $\cM$-$q$-reasonable $\gb > \ga$ with 
$j_{\gb , q} ( \gb ) = j_{\ga , q} ( \ga )$.
Suppose that $\gd \in \gD_{\ga , q}$.  
We must find $\ge \in \gD_{\gb , q}$
such that
$\gl^{\ga , q}_\gd = \gl^{\gb , q}_\ge$.
Let 
$$\jtil \ = \ j_{\gb , q}^{-1} \circ j_{\ga , q}
\ : \ \cH_{\ga , q} \lra \cH_{\gb , q}.$$
Note that $\crit (\jtil ) = \ga$ and $\jtil ( \ga ) = \gb$.
Let 
$G$ be the 
$(\ga , \gb )$-extender derived from $\jtil$.
Let
$\cP = \ult ( \cN (\gd ), G )$
and $\itil : \cN (\gd ) \lra \cP$ be the ultrapower map.
Put
$\ge = (\ga^+ )^\cP$.
We show that this $\ge$ works.

Let 
$\ktil : \cP \lra \jtil ( \cN (\gd ))$ be the natural almost
$\gS_1$-elementary embedding such that 
$\jtil \res |\cN (\gd ) | = \ktil \circ \itil$.
It is easy to see that $\ge = \crit ( \ktil ) < \gd_{\gb , q}$
and $\ktil (\ge ) = \gd_{\gb , q}$.
So $\ge$ is a local $\gb^+$ less than $\gd_{\gb , q}$,
that is,
$\ge \in \gG_\gb \cap \gd_{\gb , q}$.

The verification that $\ge \in \gD_{\gb , q}$ amounts to 
checking that various bounded subsets of 
$j_{\ga , q } (\ga ) = j_{\gb , q} (\gb )$
that are in the range of $k^{\ga , q}_\gd$ are also in the range
of $k^{\gb , q}_\ge$,
which follows from the ($\Lra$) direction of the following easy claim.

\bigskip

\noindent
{\bf Claim.}
Suppose that $\xi < j_{\ga , q } (\ga ) = j_{\gb , q} ( \gb )$.
Then
$$\exists
a \in [\gk ]^{<\go} \ \ \exists f \in J^\Evec_\gd
\ \ \ j_{\ga , q} (f ) (a  \cup s )= \xi
\ \ \ \ \Llra \ \ \ \ 
\exists
b \in [\gk ]^{<\go} \ \ \exists g \in J^\Evec_\ge
\ \ \ j_{\gb , q} (g ) (b  \cup s )= \xi$$

\bigskip

\noindent
($\Lra$) \ \ Let $b = a - \gb$ and
$c = a \cap \gb$, and $g (u ) = \jtil (f ) (c \cup u )$.
Then $j_{\gb , q} (g) ( b \cup s ) = \xi$.
Since $g \in \ran (\ktil )$ and is a subset of $\gb$, 
we have that $g \in J^\cP_\ge = J^\Evec_\ge$.

\noindent
($\Lla$) \ \ 
Again, $g \in J^\cP_\ge$, so we may
write $g = \jtil ( h ) (c )$ for some $h \in J_\gd^\Evec$
and $c \in [\gb ]^{< \go }$.
Let $a = c \cup b$ and $f (u) = h(u^{c , a \cup s} )(u^{b \cup s , a \cup s} )$.
Then 
$j_{\ga , q} (f ) (a  \cup s )= \xi$.

\bigskip

But now, using both directions of the claim, we see that $\gl^{\ga , q}_\gd
= \gl^{\gb , q}_\ge$.

\bigskip

\noindent
{\bf Clause (b).}
Suppose that $\gb$ is a limit of $\cM$-$q$-reasonable ordinals.
It is clear that clauses (i)--(iii) of Definition
1.12 hold at $\gb$.
This is enough to make sense of the definitions of 
$j_{\gb , q}$, $\gd_{\gb , q}$, and $\gl^{\gb , q}_\gd$ for 
$\gd < \gd_{\gb , q}$.
It remains to see that 1.12(iv) and 1.12(v) hold.
We show 1.12(iv),
$\gG_\gb \cap \gd_{\gb , q}$ is cofinal in $\gd_{\gb , q}$,
and leave 1.12(v) to the reader.

Since the function $\gg \mapsto j_{\gg , q } (\gg )$ is non-increasing
on the $\cM$-$q$-reasonable ordinals below $\gb$,
$\{ j_{\gg , q } (\gg ) \mid \gg < \gb$ is $\cM$-$q$-reasonable$\}$
is finite.
So, without loss of generality,
$\gb > \ga$ and
$j_{\ga , q } (\ga ) = j_{\gb , q} ( \gb )$.

Let $G$ be as in the proof of Clause (a),
and for each $\gd < \gd_{\ga , q}$, let $\cP_\gd = \ult ( \cN (\gd ), G )$
and
$\ge_\gd = (\gb^+)^{\cP_\gd}$
Let
$S = \{ \ge_\gd \mid \gd < \gd_{\ga , q} \}$.
The proof of Clause (a) shows that 
$S \subseteq \gG_\gb \cap \gd_{\gb , q}$ and that 
$\gl^{\ga , q}_\gd = \gl^{\gb , q}_{\ge_\gd}$
whenever $\gd < \gd_{\ga , q}$.
Let $\ge = \sup (S)$.
If $\ge < \gd_{\gb , q}$,
then using the proof of Lemma 1.13(d), we get 
that $\gl^{\gb , q}_\ge < \gn$.
But $\gl^{\gb , q}_\ge > \gl^{\gb , q}_{\ge_\gd}
= \gl^{\ga , q}_\gd$ for every $\ga \in \gD_{\ga , q}$,
and so $\gl^{\gb , q}_\ge \geq \gn$, which is absurd.
So 1.12(iv) holds for $\gb$.

\bigskip

\noindent
{\bf Clause (c).}
Assume that $q = p_1 (\cM )$.
Since $\cM$ is $1$-sound,
for large enough $\ga < \gk$,
$$\gk \in H^\cM_1 (\ga \cup p_1 (M)).$$
For such $\ga$, 
condition 1.12(ii) rules out $\ga$ being $\cM$-$p_1 (M)$-reasonable.

\bigskip

\noindent
{\bf Clause (d).}
That $\ga$ is $\cQ^{\ga , q}_\gd$--$p$-reasonable
follows from earlier lemmas.
Suppose, to the contrary, that $\gb$
is $\cQ^{\ga , q}_\gd$-$p$-reasonable.
It is easy to see that $i^{\ga , q}_\gd ( f )  (a) =\gk$.
Also,
$i^{\ga , q}_\gd = (j_{\ga , p})^{\cQ^{\ga , q}_\gd}$.
But then,
using the Claim in the proof of 1.17(a),
we see that $\gk \in \ran (j_{\gb , p})^{\cQ^{\ga , q}_\gd}$.
This is a contradiction,
since 
$\gb$ violates 1.12(ii) in the definition of
$\cQ^{\ga , q}_\gd$--$p$-reasonable.
\ \ $\blacksquare$ {\tiny \ \ Lemma 1.17}

\bigskip

\noindent
{\bf Definition 1.18. }
Suppose that $\gn \in \gG$ and $\cM \in \cM (\gn )$.
Define
$$e ( \cM ) 
= \{ j_{\ga , q} ( \ga ) 
\mid (\ga , q) 
\mbox{ is $\cM$-considered and } s_{\ga , q}^\cM = \emptyset \}.$$
Note that $e(\cM )$ is a finite set that
may be empty.
For each $\gz \in e ( \gn )$,
define $\ba ( \cM , \gz )$ to be the largest $\ga$ such that
$( \ga , q)$ is $\cM$-considered
and $j_{\ga , q} (\ga ) = \gz$.
Suppose that $\ga = \ba ( \cM , \gz )$.
Define
$\bb (\cM , \gz )$ to be the least $\gb \geq \ga$
such that $\gk = j_{\ga , q} (f) (a)$
for some $a \in [ \gb ]^{< \go }$
and some $f : [\ga ]^{|a|} \lra \ga$.
If $f$ is least such that 
$\gk = j_{\ga , q} (f) (a)$
for some $a \in [\bb (\cM , \gz )]^{<\go }$,
then 
define
$\gD ( \cM , \gz )$ to be the tail of $\gd \in \gD_{\ga , q}$
for which $f \in J^\Evec_\gd$.
Let 
$$S (\cM , \gz ) = 
\{ \gl^{\ga , q}_\gd \mid \gd \in \gD (\cM , \gz ) \}.$$
For any $\gd \in \gD ( \cM , \gz)$,
let
$\cQ  (\cM , \gz , \gl^{\ga , q}_\gd ) = \cQ^{\ga , q}_\gd$
and $\gZ ( \cM , \gz , \gl^{\ga , q}_\gd) = i^{\ga , q }_\gd ( \ga )$.

\bigskip

Note that 1.17(c) justifies the definition of $\ba ( \cM , \gz )$.
We remark that if $\gn \in \gG_\e$ and $\gm = \crit (E_{\gb ( \gn )})$,
then $S(\cJ_{\gb ( \gn )} , j ( \gm ))$ is a tail of $B_\e ( \gn )$.
The following lemma puts together some facts that we have already
established.  It is used, among other things, 
to justify the subsequent inductive definition.

\bigskip

\noindent
{\bf Lemma 1.19.}
Suppose that $\gn \in \gG$,
$\cM \in \cM (\gn )$, and $\gz \in e (\cM )$.
Assume that $\gnbar$ is a limit point of $S(\cM , \gz )$.
Then the following hold.
\begin{list}{}{}
\item[(a)]{$S(\cM , \gz )$ is club in $\gn$}
\item[(b)]{$\ot ( S(\cM , \gz ) )\leq \ba (\cM , \gz )^+ \leq \gk$}
\item[(c)]{$\ba ( \cM , \gz ) \leq \bb ( \cM , \gz ) < \gk$}
\item[(c)]{$\cQ ( \cM , \gz , \gnbar )
		\in \cM (\gnbar )$.}
\item[(d)]{$Z (\cM , \gz , \gnbar ) \in e (\cQ ( \cM , \gz , \gnbar ) )$}
\item[(e)]{
\begin{tabbing}
$\ba (\cM , \gz )$
\=
$\leq \ba ( \cQ ( \cM , \gz , \gnbar ), Z (\cM , \gz , \gnbar ))$\\
\>
$\leq \bb ( \cQ ( \cM , \gz , \gnbar ), Z (\cM , \gz , \gnbar ))
= \bb ( \cM , \gz )$
\end{tabbing}
}
\item[(f)]{$\gnbar \ \cap \ S(\cM , \gz ) 
\subseteq S (\cQ ( \cM , \gz , \gnbar ), Z (\cM , \gz , \gnbar ))$}
\end{list}

\bigskip

\noindent
{\bf Definition 1.20.}
We define a function $C_\e$  with domain
$$\{ (\cM , \gz ) \mid \cM \in \cM ( \gn )
\mbox{ and } \gz \in e (\cM ) \mbox{ for some } \gn \in \gG \}$$
by induction on $\gn$.
Suppose that $\gn \in \gG$.
Assume that we have defined 
$C_\e$ restricted to the set
$$\{ (\cM , \gz ) \mid 
\cM \in \cM ( \gnbar )
\mbox{ and } \gz \in e (\cM ) \mbox{ for some } \gnbar < \gn \}$$
Then for $\cM \in \cM (\gn )$ and $\gz \in e (\cM )$, we define
$C_\e (\cM , \gz )$ to be
$$S (\cM , \gz ) \ \cup \ \bigcup
\{ C_\e (\cQ (\cM , \gz , \gnbar ) , Z (\cM , \gz , \gnbar ) )
\mid
\gnbar \in \lim (S (\cM , \gz )) \}.$$

\bigskip

\noindent
{\bf Lemma 1.21.}
Suppose that $\gn \in \gG$, $\cM \in \cM (\gn )$, and $\gz \in e (\cM )$.
Then the following hold.
\begin{list}{}{}
\item[(a)]{
$\card ( C_\e ( \cM , \gz ) ) \leq \bb ( \cM , \gz )^+$}
\item[(b)]{
Suppose that $\gnbar$ is a limit point of $S(\cM , \gz )$.
Let $\cQbar = \cQ (\cM , \gz , \gnbar )$ and 
$\gzbar = Z ( \cM , \gz , \gnbar)$.
Then $C_\e ( \cQbar , \gzbar ) = \gnbar \ \cap \ C_\e (\cM , \gz )$.}
\item[(c)]{
Suppose that $\gnbar$ is a limit point of $C_\e ( \cM , \gz )$.
Then there are $\cQbar \in \cM (\gnbar )$ and $\gzbar \in e(\gnbar )$
such that
$C_\e ( \cQbar , \gzbar ) = \gnbar \cap C_\e (\cM , \gz )$.}
\end{list}

\bigskip

\noindent
{\bf Proof of 1.21.}

\noindent
{\bf Clause (a).}
By induction on $\gn$.
By the definition of $C_\e ( \cM , \gz )$, 
Lemma 1.19, and the induction hypothesis,
$C_\e ( \cM , \gz ) $ is the union of 
at most $\bb (\cM , \gz )^+$ many sets, each of size 
at most $\bb (\cM , \gz )^+$.
Recalling that $\gk$ is a limit cardinal, (a) follows.

\bigskip

\noindent
{\bf Clause (b).}
By induction on $\gn$.  Assume that clause (b) holds for all $\gn' < \gn$.
Suppose that $\gnbar$ is a limit point of $S(\cM, \gz )$.

If $\gnbar$ is the greatest limit point of $S (\cM , \gz )$,
then by the induction hypothesis and the definition of $C_\e (\cM , \gz )$,
$$ C_\e ( \cM , \gz ) =
C_\e ( \cQbar , \gzbar ) \cup (S(\cM , \gz) - \gnbar),$$
and it is clear that clause (b) holds.

Consider an arbitrary $\gn' > \gnbar$ such that $\gn'$ is also a 
limit point of $S(\cM, \gz )$.
Let $\cQ' = \cQ ( \cM, \gz , \gn' )$ and $\gz' = Z ( \cM , \gz , \gn' )$.
By Lemma 1.19,
$\gn' \cap S (\cM , \gz ) \subseteq 
S(\cQ' , \gz' )$.
Thus $\gnbar$ is a limit point of 
$S (\cQ' , \gz' )$.
By the induction hypothesis, clause (b) holds at $\gn'$, and hence
$$(\star) \ \ \ \ \ C_\e (\cQbar , \gzbar )
= \gnbar \ \cap \ C_\e (\cQ' , \gz').$$

If $\gn$ is a limit of limit points of $S(\cM, \gz )$,
then clause (b) holds by ($\star$) and the definition
of $C_\e (\cM , \gz )$.

If $\gn$ is not a limit of limit points of $S(\cM, \gz )$,
then take $\gn'$ to be the greatest limit point of $S(\cM, \gz )$.
Applying ($\star$) to the definition of $C_\e (\cM , \gz )$ again,
we see that clause (b) holds.

\bigskip
\noindent
{\bf Clause (c).}
By induction on $\gn$.
If $\gnbar$ is a limit point of $S (\cM , \gz )$,
then we are done by clause (b).
So we may assume that there is a $\gn' > \gn$
such that $\gn'$ is a limit point of $S (\cM , \gz )$.
Let 
$\cQ' = \cQ ( \cM , \gz , \gn' )$ and $\gz' = Z ( \cM , \gz , \gn' )$.
Then
$C_\e ( \cQ' , \gz' ) = \gn' \cap C ( \cM , \gz )$
by clause (b).
So $\gnbar$ is a limit point of $C_\e ( \cQ' , \gz' )$.
By the induction hypothesis,
there is an $\cMbar \in \cM (\gnbar )$ and $\gzbar \in e(\cMbar)$
such that 
$$C_\e ( \cQbar , \gzbar ) = \gnbar \cap C_\e (\cQ' , \gz' )
= \gnbar \cap C_\e ( \cM , \gz ),$$
so we are done. \ \ $\blacksquare$ {\tiny \ \ Lemma 1.21}

\bigskip

In the next few lemmas and definitions
we are mainly concerned with the case 
$s_{\ga , q} \not= \emptyset$.

\bigskip

\noindent
{\bf Definition 1.22.}
Suppose that $\gn \in \gG$ and $\cM \in \cM (\gn )$.
Suppose that $(\ga , q)$ is $\cM$-considered and that 
$\gd \in \gD_{\ga , q}$.  We say that $\gd$ is {\bf $F_{\ga , q}$-closed}
if and only if for every 
$x \subseteq \ga$, if there are $a \in [ \gd ]^{< \go}$ and 
$f \in J^\Evec_\gd$ such that 
$x = j_{\ga , q} ( f ) ( a \cup s_{\ga , q} )$,
then $x \in J^\Evec_\gd$.

\bigskip

Note that Definition 1.22 is consistent with Definition 1.5.
The following lemma is proved as was Lemma 1.6.

\bigskip

\noindent
{\bf Lemma 1.23.}
Suppose that $\gn \in \gG$ and $\cM \in \cM (\gn )$.
Suppose that $(\ga , q)$ is $\cM$-considered and that 
$\gd \in \gD_{\ga , q}$ is $F_{\ga ,q}$-closed.
Then $r_\ga (\cQ^{\ga , q}_\gd ) = i^{\ga , q}_\gd ( p_1 ( \cN ( \gd )))$.

\bigskip

If $\gn \in \gG_\n$ with $\gm = \crit (E_{\gb ( \gn )})$,
then we defined 
$\Phi_{\gm , \emptyset}$ to be the tail of the $F_{\gm , \emptyset}$-closed
$\gd$ in $\gD_{\gm , \emptyset}$ such that
$\gm$ is switching for $\cQ^{\gm , \emptyset}_\gd$.
(This is Definition 1.9 in the ``new notation''.) \ 
And we defined $B_\n ( \gn ) = \{ \gl^{\gm , \emptyset}_\gd \mid
\gd \in \Phi_{\gm , \emptyset} \}$.
For consistency with the definition in the next paragraph,
let us, in this case,
put $\Phi ( \cJ_{\gb ( \gn )} , \gm ) = \Phi_{\gm , \emptyset}$
and 
$C_\n ( \cJ_{\gb ( \gn )} , \gm ) = B_\n ( \gn )$.

Now suppose that $\gg \in \gG$ and that $\cM \in \cM ( \gn )$.
Suppose that $\ga \in \mbox{SC}$.
Put $q = r_\ga ( \cM )$.
Suppose that $\ga$ is $\cM$-$q$-reasonable,
and that $\ga$ is switching for $\cM$.
In particular, $s_{\ga , q} \not= \emptyset$.
Suppose that 
$A = \{ \gd \in \gD_{\ga , q} \mid \gd \mbox{ is $F_{\ga , q}$-closed}\}$
is club in $\gd_{\ga , q}$.
Suppose that for sufficiently large
$\gd$ in $A$,
we have that
$\ga$ is switching for $\cQ_\gd^{\ga , q}$.
In that case, 
we let $\Phi ( \cM , \ga )$ be the tail of $\gd$ in $A$
for which $\ga$ is switching for $\cQ_\gd^{\ga , q}$,
and we put
$$C_\n ( \cM , \ga ) 
= \{ \gl^{\ga , q}_\gd \mid \gd \in \Phi ( \cM , \ga ) \}.$$

Let us say that $\ga$ is {\bf $\cM$-critical} if and only if
$C_\n (\cM , \ga )$ is defined (that is, if it is defined by one of the
two cases above).  There are only finitely many $\cM$-critical $\ga$ 
for a given
$\cM$,
since to each critical $\ga$,
we can associate a distinct partition of $p_1 (\cM )$.
The following is an easy consequence of Lemmas 1.16(b), 
1.23, and the definitions
just given.

\bigskip

\noindent
{\bf Lemma 1.24.}
Suppose that $\gn \in \gG$ and $\cM \in \cM (\gn )$.
Suppose that $\ga$ is $\cM$-critical and that $\gnbar$ is a limit point of 
$C_\n ( \cM , \ga )$.  Say $\gnbar = \gl^{\ga , q}_\gd$, where
$q = \emptyset$ if $\ga = \gmdot^\cM$, and $q = r_\ga ( \cM )$ otherwise.
Then $\ga$ is $\cQ^{\ga , q}_\gd$-critical and
$$C_\n ( \cQ^{\ga , q}_\gd , \ga ) = \gnbar \cap C_\n (\cM , \ga ).$$

\bigskip

We are finally ready to write down our $\square_{\gk , < \go} (\gG )$
sequence, $\cF$.
If $\gn \in \gG_\h$, then
\begin{tabbing}
$\cF (\gn )$ $=$
$\{ B_\h ( \gn ) \}$ \= $\cup$ 
$\{ C_\e ( \cM , \gz ) \mid \cM \in \cM (\gn ) \ \land 
\ \gz \in e( \cM ) \}$\\
\> $\cup$ 
$\{ C_\n ( \cM , \ga ) \mid
\cM \in \cM (\gn ) \ \land \ \ga \mbox{ is $\cM$-critical}\}$
\end{tabbing}
If $\gn \in \gG_\t$, then
\begin{tabbing}
$\cF (\gn )$ $=$
$\{ B_\t ( \gn ) \}$ \= $\cup$ 
$\{ C_\e ( \cM , \gz ) \mid \cM \in \cM (\gn ) \ \land 
\ \gz \in e( \cM ) \}$\\
\> $\cup$ 
$\{ C_\n ( \cM , \ga ) \mid
\cM \in \cM (\gn ) \ \land \ \ga \mbox{ is $\cM$-critical}\}$
\end{tabbing}
If $\gn \in \gG_\l$, then
\begin{tabbing}
$\cF (\gn )$ $=$
\= $\{ C_\e ( \cM , \gz ) \mid \cM \in \cM (\gn ) \ \land
\ \gz \in e( \cM )\}$\\
\> $\cup$
$\{ C_\n ( \cM , \ga ) \mid
\cM \in \cM (\gn ) \ \land \ \ga \mbox{ is $\cM$-critical}\}$
\end{tabbing}

That $\cF$ witnesses $\square_{\gk , < \go} (\gG )$ now follows 
easily from our various lemmas and observations;
we emphasize just a few points.
First, we claim that each $\cF (\gn )$ is non-empty.  
This is visibly so if $\gn \in \gG_\h \cup \gG_\t$,
so suppose that $\gn \in \gG_\l$.  Let $\gm =\crit (E_{\gb ( \gn )})$
and $j : \LofE \lra \ult ( \LofE , E_{\gb ( \gn )} )$.
If $\gg \in \gG_\n$, then
$\gm$ is $\cJ_{\gb (\gn )}$-critical,
so $C_\n ( \cJ_{\gb ( \gn )} , \gm )$ is defined.
If $\gg \in \gG_\e$,
then $j_{\gm , \emptyset} ( \gm ) \in e(\cJ_{\gb (\gn )} )$
and $C_\e ( \cJ_{\gb ( \gn )} , j ( \gm ) )$ is defined.
Second, we claim that $\cF (\gn ) $ is finite.
This is because $\cM (\gn )$ is finite,
and $e(\cM )$ is finite whenever $\cM \in \cM (\gn )$,
and there are only finitely many $\cM$-critical ordinals.
Third, each element of $\cF (\gn )$ is a club subset of $\gn$ of order type
at most $\gk$.
This is essentially a consequence of Lemma 1.16 and Lemma 1.21(a).
Finally, the coherence property, Definition 0(d), holds
for sets $C \in \cF (\gn )$.
If $C \in \ran ( B_\h)$, then this follows from \S2 of \cite{Sch}.
If $C\in \ran ( B_\t )$, then 
there is nothing to check, since $\ot (C ) = \go$.
If $C\in \ran (C_\e )$, then we quote Lemma 1.21(c).
And if $C\in \ran (C_\n )$, then
use Lemma 1.24.
\ \ $\blacksquare$ {\tiny \ \ Theorem 1}


\bigskip

\begin{thebibliography}{MiSchSt34}

\bibitem[BM]{BM}S. Ben-David and M. Magidor,
                {\em The Weak ${\square}^*$ is Really Weaker 
                Than the Full $\square$},
                J. Symbolic Logic.  {\bf 51} 4 (1986) 1029-1033.

\bibitem[J1]{J1}R.B. Jensen,
              {\em The Fine Structure of the Constructible Hierarchy},
              Annals of Math. Logic {\bf 4} (1972) 229-308.

\bibitem[J2]{J2}R.B. Jensen,
		{\em Some remarks on $\square$ below $0^{\P}$},
		handwritten.

\bibitem[MiSchSt]{MiSchSt}W.J. Mitchell, E. Schimmerling, and J.R. Steel,
                {\em The Covering Lemma up to One Woodin Cardinal},
		  submitted for publication.

\bibitem[MiSt]{MiSt}W.J. Mitchell and J.R. Steel,
                  {\em Fine Structure and Iteration Trees},
		  Lecture Notes in Logic 3, Springer-Verlag (1994).

\bibitem[Sch]{Sch}E. Schimmerling,
		{\em Combinatorial Principles in the Core Model for one
		Woodin Cardinal},
		to appear in Annals of Pure and Applied Logic.

\bibitem[MiSch]{MiSch}W.J. Mitchell and E. Schimmerling,
		{\em Weak Covering Without Countable Closure}.

\bibitem[SchSt]{SchSt}E. Schimmerling and J.R. Steel,
		{\em Fine Structure for Tame Inner Models},
		submitted for publication.

\bibitem[S]{S}R.M. Solovay,
		{\it The fine structure of L[$\mu$]},
		unpublished.

\bibitem[SRK]{SRK}R.M. Solovay, W.N. Reinhart, and A. Kanamori,
		{\em Strong axioms of infinity and elementary embeddings},
		Ann. Math. Logic, {\bf 13} (1978), 73-116.

\bibitem[St1]{St1}J.R. Steel,
              {\em The Core Model Iterability Problem},
	      to appear in Lecture Notes in Logic.

\bibitem[St2]{St2}J.R. Steel,
	{\em Inner Models With Many Woodin Cardinals},
	Ann. of Pure and Appl. Logic 65 (1993), 185--209.

\bibitem[St3]{St3}J.R. Steel,
	{\em Core Models With More Woodin Cardinals},
	in preparation.

\bibitem[T]{T}S. Todorcevic,
              {\em A Note on the Proper Forcing Axiom},
              Contemporary Mathematics {\bf 95} (1984) 209-218.

\bibitem[We]{We}P.D. Welch,
                {\em Combinatorial Principles in the Core Model},
                Oxford doctoral dissertation.

\bibitem[Wy]{Wy}D.J. Wylie,
                MIT doctoral dissertation (1990).

\end{thebibliography}
\end{document}